\documentclass[11pt]{article}
\usepackage{amssymb,latexsym}
\usepackage{graphicx}

\newtheorem{theorem}{Theorem}[section]

\newtheorem{defi}{Definition}[section]

\def\slfrac#1#2{\hbox{\kern.1em %
 \raise.5ex\hbox{\the\scriptfont0 #1}\kern-.11em %
 /\kern-.15em\lower.25ex\hbox{\the\scriptfont0 #2}}}

\newcommand{\eqn}[1]{(\ref{#1})}

\newcommand{\eeq}{\end{equation}}
\newcommand{\beql}[1]{\begin{equation}\label{#1}}

\newcommand{\ZZ}{{\mathbb Z}}
\newcommand{\RR}{{\mathbb R}}

\newcommand{\CC}{{\mathbb C}}

\newcommand{\sF}{{\cal F}}
\newcommand{\sG}{{\cal G}}

\newcommand{\dM}{{N}^{1}}
\newcommand{\dN}{{D}}
\newcommand{\w}{{y}}

\newcommand{\tauS}{\tau^{S}}
\newcommand{\tb}{\bar{b}}

\newcommand{\ddL}{\Omega^{L}}

\newcommand{\muT}{{\mu_S}}
\renewcommand{\ll}{\langle\!\langle}
\renewcommand{\gg}{\rangle\!\rangle}

\makeatletter
\def\@sect#1#2#3#4#5#6[#7]#8{\ifnum #2>\c@secnumdepth
     \def\@svsec{}\else
     \refstepcounter{#1}\edef\@svsec{\csname the#1\endcsname.\hskip .75em }\fi
     \@tempskipa #5\relax
      \ifdim \@tempskipa>\z@
        \begingroup #6\relax
          \@hangfrom{\hskip #3\relax\@svsec}{\interlinepenalty \@M #8\par}%
        \endgroup
       \csname #1mark\endcsname{#7}\addcontentsline
         {toc}{#1}{\ifnum #2>\c@secnumdepth \else
                      \protect\numberline{\csname the#1\endcsname}\fi
                    #7}\else
        \def\@svsechd{#6\hskip #3\@svsec #8\csname #1mark\endcsname
                      {#7}\addcontentsline
                           {toc}{#1}{\ifnum #2>\c@secnumdepth \else
                             \protect\numberline{\csname the#1\endcsname}\fi
                       #7}}\fi
     \@xsect{#5}}
\def\@begintheorem#1#2{\it \trivlist \item[\hskip \labelsep{\bf #1\ #2.}]}

\def\plain{plain}\ifx\fmtname\plain\csname fi\endcsname
     
     \input docstrip
     \preamble

     Do not distribute the stripped version of this file.
     The checksum in the header refers to the documented version.

     \endpreamble
     \generateFile{here.sty}{t}{\from{here.doc}{}}
     \endinput
\fi
\ifcat a\noexpand @\let\next\relax\else\def\next{%
    \documentstyle[here,doc]{article}\MakePercentIgnore}\fi\next
\ifx\@Hxfloat\@Hundef\else\expandafter\endinput\fi
\let\@Hxfloat\@xfloat
\def\@xfloat#1[{\@ifnextchar{H}{\@HHfloat{#1}[}{\@Hxfloat{#1}[}}
\def\@HHfloat#1[H]{%
\expandafter\let\csname end#1\endcsname\end@Hfloat
\vskip\intextsep\vbox\bgroup\def\@captype{#1}\parindent\z@
\ignorespaces}
\def\end@Hfloat{\egroup\vskip \intextsep}

\makeatother

\catcode`\@=11
\renewcommand{\section}{
        \setcounter{equation}{0}
        \@startsection {section}{1}{\z@}{-3.5ex plus -1ex minus
        -.2ex}{2.3ex plus .2ex}{\large\bf}
        }
\catcode`@=12

\begin{document}

\begin{center}
{\Large 
{\bf    The Takagi Function and Its Properties} 
}\\

\vspace{1.5\baselineskip}
{\em Jeffrey C. Lagarias
\footnote{This work was supported by NSF Grants DMS-0801029 and DMS-1101373.} }\\
\vspace*{.2\baselineskip}
Dept. of Mathematics \\
University of Michigan \\
Ann Arbor, MI 48109-1043\\
\vspace*{1.5\baselineskip}

\vspace*{2\baselineskip}
(January 2,  2012) \\
\vspace{3\baselineskip}
{\bf ABSTRACT}
\end{center}
The Takagi function $\tau(x)$ is a continuous non-differentiable function
introduced by Teiji Takagi  in 1903. 
It has appeared in a surprising number of different mathematical 
contexts, including mathematical analysis, probability theory and number theory.
 This paper surveys properties of this function.

%
%
%
%
\setlength{\baselineskip}{1.0\baselineskip}

\section{Introduction}

\begin{figure}[h]
\centering
\includegraphics[height=2.5 in]{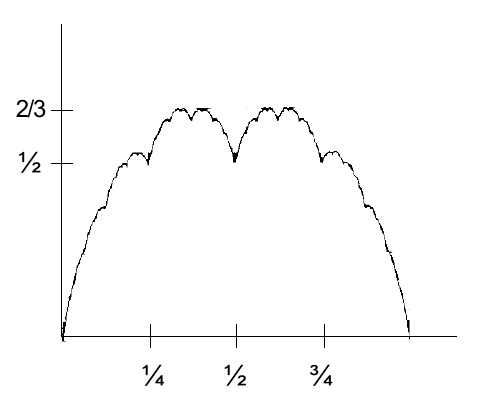}
\caption{Graph of the Takagi function $\tau(x)$.}
\label{fig11}
\end{figure}

The Takagi function $\tau(x)$ was introduced by T.Takagi \cite{Tak03}
in 1903 as an example of an everywhere non-differentiable function
on $[0,1]$. It can be defined on the unit interval $ x\in [0,1]$
by
\begin{equation}\label{eq101}
\tau(x) := \sum_{n=0}^{\infty} \frac{1}{2^n}\ll 2^n x \gg
\end{equation}
where $\ll x \gg$ is the distance from $x$ to the nearest integer. 
Takagi  defined it using binary expansions, 
and showed that his
definition was consistent for numbers having two binary expansions
(dyadic rationals). The function is pictured in Figure \ref{fig11}.

%
%

An  immediate generalization of the Takagi function  is to  set,
 for integer $r \ge 2$, 
\begin{equation}\label{eq102a}
\tau_r(x) :=\sum_{n=0}^{\infty} \frac{1}{r^n} \ll r^n x \gg.
\end{equation}
In 1930 van der Waerden \cite{VDW30} studied the 
function $\tau_{10}$ and proved its non-differentiability.
In 1933 Hildebrandt \cite{Hi33} simplified his construction
to rediscover the Takagi function.
Another rediscovery of the Takagi function was made by
of de Rham \cite{deR57} in 1957. 

The Takagi function has appeared in a surprising number of different
contexts, including analysis, probability theory and number theory. 
It is a prescient example of a self-similar construction.
An important feature of this function 
is that it satisfies many self-similar functional equations,
 some of which are related to the dilation equations that
appear in wavelet theory.
It serves as a kind of exactly solvable 
``toy model" representing a solution to a discrete version of the Laplacian operator.
It serves as an interesting test case for 
determining various measures of irregularity of behavior of a function.
It also turns out to be related to a number of interesting singular functions.

This object of this paper is to survey properties of Takagi function
across all these fields.
In Section \ref{sec2} we briefly review
some  history of work on the Takagi function. 
  In Sections \ref{sec3} and \ref{sec4} 
we review basic analytic properties of the Takagi function.
In Section \ref{sec5} we discuss  its connection to dynamical systems.
In Section \ref{sec6} we treat its Fourier transform, and
 in Section \ref{sec7} its relation with  Bernoulli convolutions
in probability theory.
Section \ref{sec8}  summarizes analytic results quantifying the
local oscillatory behavior of the Takagi function
at different scales, which imply
its non-differentiabilty.
Results  related to number theory appear 
in Sections \ref{sec9} and  \ref{sec10}, connecting
it with binary digit sums and with the Riemann hypothesis,
respectively.
In Section \ref{sec11}  we describe 
 properties of the graph of the Takagi function.
In Sections \ref{sec12} - \ref{sec14} we present  results 
concerning the 
level sets of the Takagi function.  These include results  
 recently obtained jointly with 
Z. Maddock (\cite{LM10a}, \cite{LM10b}).
Many of the results in these sections are based on
properties specific to the Takagi function. This  again exhibits
its value as a ``toy model", where exact calculations are possible. 

In this paper most  results are given without proof. We have included proofs
of several  results  not conveniently available, and some of these  may
not have been noted before, e.g. Theorem \ref{th32a} and Theorem
\ref{th82}.
We also draw the reader's attention to another recent survey 
of results on the Takagi function given by  Allaart and Kawamura \cite{AK11},
which is somewhat complementary to this one.

%
%
%
%

\section{History}\label{sec2}

The Takagi function was introduced during a period when
the general structure of non-differentiable functions was
being actively explored. This started  from 
Weierstrass's  discovery of an everywhere
non-differentiable function, which he lectured on 
as early as 1861, but which was first published
(and attributed to Weierstrass)
 by Du Bois-Reymond \cite{dBR1875}
in 1875. Weierstrass's example has no finite
or infinite derivative at any point. 
The  example of Takagi is much simpler,
and has no point with a finite
derivative,  but it does have some points with a well-defined
infinite derivative, see Theorem ~\ref{th86a} below.
Pinkus \cite{Pin00} gives additional history on this problem.\smallskip

In 1918  Knopp \cite{Knopp18} studied many variants of
non-differentiable functions and reviewed earlier work, 
including that of
Faber \cite{Fab08}, \cite{Fab10a}  and Hardy \cite{Har16},
among others.
He considered functions of the general form
\begin{equation}\label{eq131}
F(x) := \sum_{n=0}^{\infty} a^n \, \phi( b^n x),
\end{equation}
where $\phi(x)$ is a given periodic  continuous function of period one,
for real numbers $0< a< 1$ and $b>1$.
This general form 
includes the Takagi function as well as functions in the  Weierstrass 
non-differentiable function family
$$
W_{a, b}(x) := \sum_{n=0}^{\infty} a^n \cos (\pi b^n x).
$$
Weierstrass showed this function has no finite or infinite derivative
when  $0<a<1$ and $b$ is an odd integer with
 $ab >1+ \frac{3 \pi}{2}.$
In 1916 Hardy \cite[Theorem 1.31] {Har16} proved that the Weierstrass
 function has no finite derivative for real
$a, b$  satisfying $0 < a<1$ and $ab \ge 1$.  For the special case where
$\phi(x) = \ll x\gg$, relevant to the Takagi function,
 Knopp proved (\cite[p. 18]{Knopp18}) the non-differentiability 
for real  $0< a<1$ and $b$ being an even integer 
with   $ab>4$.  The Takagi function has $ab=1$, so is not covered by
Knopp's result; however 
 a later result of  Behrend \cite[Theorem III]{Beh48} in 1948, 
applies to  $\phi(x) = \ll x\gg$ 
and establishes  no finite derivative for  integer $b$
and $0< a<1$,  having  $ab\ge1$,
with some  specific exceptions.

In the 1930's there were significant developments in probability
theory, including  its formalization in terms of measure
theory. The expression of Lebesgue measure $[0, 1]$ in terms of the
induced measure on the binary coefficients, reveals
that these measures are given by independent coin flips (Bernoulli trials). 
This measure density can be expressed as 
an infinite product (Bernoulli convolution), see Kac \cite{Kac59}
for a nice treatment.
In 1934,  Lominicki and Ulam \cite{LU34} 
studied related measures where biased coin flips are allowed. These  measures
are generally singular with respect to Lebesgue measure.
In 1984 Hata and Yamaguti \cite{HY84} noted a relation 
of the Takagi function to this family of measures, stated  below
in Theorem~\ref{th22a}. 

There has been much further work studying properties of  Bernoulli convolutions, 
as well as more general
infinite convolutions, and their associated measures.
 Additional  motivation comes from work of 
Jessen and Wintner \cite{JW35} concerning the Riemann zeta function, 
which  is  described at length in
1938 lecture notes of Wintner \cite{Wintner38}.
 Peres, Schlag and Solomyak  \cite{PSS00} give a recent progress report on
Bernoulli convolutions, and Hilberdink \cite{Hi04} surveys connections of
Bernoulli convolutions with
analytic number theory.

In the 1950's Georges de Rham (\cite{deR56}, \cite{deR56b}, \cite{deR57})
considered self-similar constructions of geometric objects,
again constructing a function equivalent to the Takagi function.
In this context Kahane \cite{Kah59} noted an important property of the
Takagi function. Similar constructions appear in the theory of splines, 
of functions iteratively
constructed using control points. 

The Takagi function appeared in number theory in connection with
 the summatory functions of various arithmetic
functions associated to binary digits.  The analysis of such
sums began with Mirsky \cite{Mir49} in 1949, but  the connection 
with the Takagi function was 
first observed by Trollope \cite{Trollope68} in 1968. It was further
explained in a very influential paper of Delange \cite{Delange75}
in 1975. A reformulation of the theory in terms of Mellin transforms
was given by Flajolet et al \cite{FGKPT94}. These results are discussed
in Section \ref{sec8}.

The Takagi function can also be viewed in terms of a dynamical
system, involving iterations of  the tent map. This  viewpoint was taken
 in 1984 by Hata and Yamaguti \cite{HY84}.
Here the Takagi function  can be  seen as a kind of fractal.
For further information 
on the fractal interpretation see Yamaguti, Hata and Kigami \cite{YHK97}.

In the 1990's 
the construction of compactly-supported wavelets  led to
the study of {\em dilation equations}, which are functional 
equations which linearly relate  functions at two (or more) different scales,
see Daubechies \cite{Dau92}.
Basic results on the solution of such
equations appear in Daubechies and the author \cite{DL91}, \cite{DL92a}.
These functions can be described in terms of infinite products, which
are generalizations of Bernoulli convolutions (see \cite{DL92b}). 
The Takagi function appears in this general context
because it satisfies a non-homogenous  dilation equation,
driven by an auxiliary function, which is  stated in Theorem ~\ref{th22}.
Related connections with de Rham's functions were observed
by Berg and Kr\"{u}ppel \cite{BK00} in 2000 and 
 by  Kr\"{u}ppel \cite{Kr09} in 2009.
\smallskip

The Takagi function has appeared in additional contexts. 
In 1995 Frankl, Matsumoto, Rusza and Tokushige \cite{FMRT95} gave
a combinatorial application. For a family $\sF$ of $k$-element sets of the
$N$-element set $[N] := \{ 1, 2, ..., N\}$ and for  any $\ell< k$
 the {\em shadow} $\Delta_{\ell}(\sF)$ of $\sF$ on $\ell$-element sets is the set 
of all $\ell$-element sets that are contained in some set in $\sF$.
The Kruskal-Katona theorem asserts that the minimal size of $\ell$-shadows
of  size $m$ families $\sF$ of $k$-element sets is attained  by choosing the  sets of $\sF$
as the first $m$ elements in the  $k$-element subsets of $[N]$ 
ordered in  the co-lexicographic order. 
The {\em Kruskal-Katona function} gives this number, namely
$$
K_{l}^k(m) := \min \{ \#\Delta_{\ell}(\sF) : \sF ~\mbox{consists of  k-element sets},~~|\sF|= m\}.
$$
This number is independent of the value of $N$, requiring only that $N$ be sufficiently
large that such  families exist, i.e. that ${N\choose k} \ge m$.
The {\em shadow function} $S_k(x)$ is a normalized version of the  Kruskal-Katona function, taking $\ell = k-1$,
which is given by
$$
S_k(x) := \frac{k}{{{2k-1}\choose{k}}}
 K_{k-1}^k( \lfloor {{2k-1}\choose{k}} x\rfloor), ~~~0 \le x \le 1.
$$
Theorem 4 of  \cite{FMRT95} states  that as $k \to \infty$ the shadow functions $S_k(x)$ uniformly converge
to the Takagi function $\tau(x)$.
The Takagi  function also appears in models of diffusion-reaction processes (\cite{GK98})
and in basins of attraction of dynamical systems (\cite{YTM88}).

%
%
%
%

\section{Basic Properties: Binary Expansions} \label{sec3}

Takagi's definition of his function was
in terms of binary expansions, which we write
\beql{221}
x = \sum_{j=1}^{\infty} \frac{b_j}{2^j}=0.b_1b_2b_3 \cdots, ~~~~ \mbox{each}~b_j \in \{0, 1\}. 
\eeq
The binary expansion of $x$ is unique except for 
dyadic rationals $x= \frac{k}{2^n}$, which have two possible expansions.
For $0 \le x \le 1$ the distance to the nearest integer function $\ll x \gg$ is
\beql{204a}
\ll x \gg  ~:= \left\{ 
\begin{array}{lcl}
x & \mbox{if}& 0 \le x  < \frac{1}{2},~~\mbox{i.e.} ~b_1=0\\
1-x & \mbox{if}&  \frac{1}{2} \le x \le 1, ~~\mbox{i.e.} ~ b_{1}=1.
\end{array}
\right.
\eeq
For  $n \ge 0$, we have 
\beql{204b}
\ll 2^n x \gg ~= \left\{ 
\begin{array}{lcl}
0.b_{n+1} b_{n+2} b_{n+3} ... & \mbox{if}& b_{n+1}=0\\
~&~&~\\
0.\tb_{n+1} \tb_{n+2} \tb_{n+3} ...& \mbox{if}&   b_{n+1}=1,
\end{array}
\right.
\eeq
where we use the bar-notation
\beql{204c}
\bar{b}= 1-b , ~~~\mbox{for}~~~ b=0 ~\mbox{or} ~1,
\eeq
to mean complementing a bit.

%
\begin{defi}\label{de23} 
{\em 
Let $x \in [0,1]$ have binary expansion
$x= \sum_{j=1}^{\infty} \frac{b_j}{2^j}= 0.b_1 b_2 b_3...$, with each $b_j \in \{0, 1\}$. 
For each $j \ge 1$ we define the following integer-valued functions.

(1) The {\em digit sum function} $\dM_j(x)$ is
\beql{221a}
\dM_j(x) := b_1 +b_2 + \cdots + b_j.
\eeq
 We also let $N_j^{0}(x)  = j - \dM_j(x)$ count the number of $0$'s in the 
 first $j$ binary digits of $x$.
 
 \smallskip
(2) The {\em  deficient digit function} $\dN_j(x)$ is given by
\beql{222}
\dN_j(x):= N_j^{0}(x) - \dM_j(x) = j- 2\dM_j(x)  = j-  2(b_1+b_2+ \cdots + b_j) .
\eeq
Here we use  the convention that $x$ denotes a binary expansion; 
dyadic rationals have two different binary expansions, and all functions 
$N_j^0(x)$, $\dM_j(x)$, $\dN_j(x)$ depend
on which binary expansion is used.
(The name ``deficient digit function" reflects the
fact that  $\dN_j(x)$ counts the excess of binary digits $b_k=0$ over those with $b_k=1$
in the first $j$ digits, i.e. it is positive if there are more $0$'s than $1$'s.  )
}
\end{defi}

Takagi's original characterization of his function, which he used
to establish non-differentiability,  is as follows.

%
%
%
%

\begin{theorem}~\label{le21} {\rm (Takagi 1903)}
For $x= 0.b_1b_2 b_3 ...$ the Takagi function is given by
\beql{201}
\tau(x) = \sum_{m=1}^{\infty} \frac{\ell_m}{2^m},
\eeq
in which $0 \le \ell_m= \ell_m(x) \le m-1$ is the integer
\beql{202}
\ell_m(x) = \# \{ i:~ 1 \le i < m,~~b_i \ne b_{m} \}.
\eeq
In terms of the digit sum function $\dM_m(x)= b_1 + b_2 + ...+b_m$,
\beql{203}
\ell_{m+1} (x) = \left\{ \begin{array}{lcl}
\dM_{m}(x) & \mbox{if}& b_{m+1}=0,\\
~&~&~\\
m-\dM_{m}(x) & \mbox{if}& b_{m+1}=1.
\end{array}
\right.
\eeq
\end{theorem}

Dyadic rationals $x= \frac{k}{2^m}$ have two binary expansions,
and in consequence the formulas above give two expansions for
$\tau(x)$.  Theorem \ref{le21} asserts that these
expansions give the same value; one may verify that $\tau(x)$ itself will then be another dyadic
rational, with the same or smaller denominator. See \cite[Lemma 2.1]{LM10a} for a proof
of this result.
\smallskip

We  deduce some basic properties of the Takagi function from its
definition.

%
\begin{theorem} ~\label{th28}
(1) The Takagi function $\tau(x)$ maps rational numbers $x$ to
rational numbers $\tau(x).$

(2) The values of the Takagi function satisfy 
$0 \le \tau(x) \le \frac{2}{3}.$  The minimal value  $y=0$ is attained only at
$x=0, 1.$ The maximal value $y= \frac{2}{3}$  is also attained
at some rational $x$, in particular  $\tau(\frac{1}{3}) = \frac{2}{3}$. 
 \end{theorem}

\paragraph{Proof.} (1) This follows from  \eqn{eq101}
because for rational $x$ the sequence $\ll 2^nx \gg$ 
takes rational values with bounded denominators, and becomes eventually periodic.
Summing geometric series then gives the rationality.

(2) The lower bound case is clear, and by inspection is attained for $x=0,1$
only.  The upper bound is proved by checking that 
$\ll x \gg +\frac{1}{2} \ll 2x \gg \le \frac{1}{2}$ holds  for all $x \in [0, 1]$,
and using this on successive pairs of terms in \eqn{eq101} to get 
$\tau(x) \le \frac{1}{2} + \frac{1}{8} + \frac{1}{32} + \cdots = \frac{2}{3}.$
One checks from \eqn{eq101} that $\tau(\frac{1}{3}) = \frac{2}{3}$. 
$~~~\Box$\medskip

\paragraph{Remark.}
The set of values $x$ having $\tau(x)= \frac{2}{3}$ is quite large; see Theorem \ref{th92}.
Concerning the  converse direction to  Theorem~\ref{th28} (1): It is not known which rationals
$y$ with $0 \le y \le \frac{2}{3}$ have the property that
there is some rational $x$ such that $\tau(x) =y$. \bigskip

\begin{figure}[h]
  \begin{center}$                                                               
    \begin{array}{ccc}
      \includegraphics[width=1.75in]{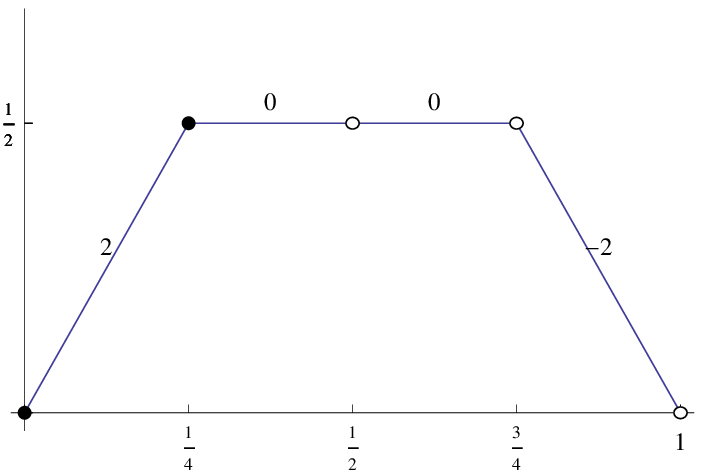} &
      \includegraphics[width=1.75in]{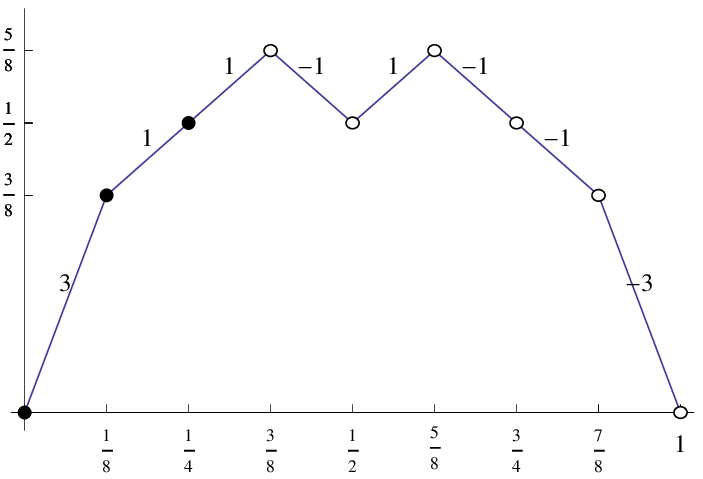} &
      \includegraphics[width=1.75in]{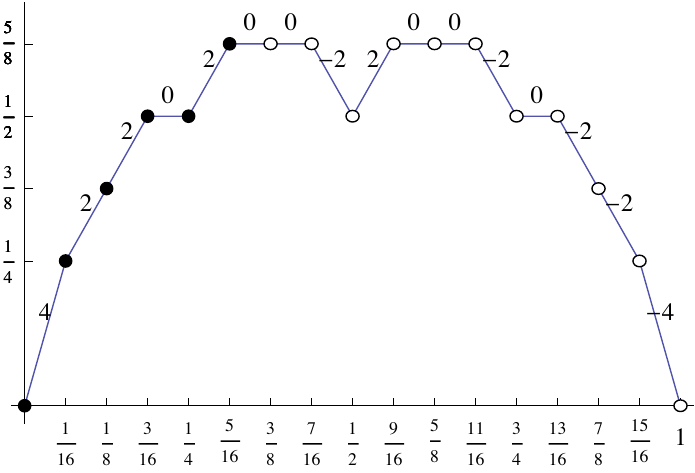} 
     \end{array}$
  \end{center}
  \caption{Approximants to Takagi function:
  (left to right) $ \tau_2(x), \tau_3(x), \tau_4(x)$. 
      Slopes of linear segments are labelled on graphs.
    A vertex $(x,y)$ is 
    marked by a solid point if and only if $x \in \Omega^L$,
    as defined in \eqn{1221}. }
  \label{fig21}
\end{figure}

The Takagi function $\tau(x)$ can be constructed 
as a limit of piecewise linear approximations. 
The  {\em partial Takagi function} of level $n$ is given by:
\begin{equation}\label{220a}
\tau_n(x) := \sum_{j=0}^{n-1}  \frac{ \ll 2^j x\gg}{2^j},
\end{equation}
See Figure \ref{fig21} for $\tau_2(x), \tau_3(x)$ and $\tau_4(x)$.

%
%
%
%

\begin{theorem}~\label{le23}
The piecewise linear function
$\tau_n(x) = \sum_{j=0}^{n-1}  \frac{ \ll 2^j x\gg}{2^j}$
 is  linear on each dyadic
interval $[\frac{k}{2^n}, \frac{k+1}{2^n}]$.

(1) On each such interval $\tau_n(x)$ has integer slope between $-n$ and $n$ given by
the deficient digit function  
$$
\dN_n(x) = N_{n}^0(x) - \dM_n(x) =n - 2(b_1+b_2 + \cdots + b_n),
$$
Here  $x=0.b_1b_2 b_3...$ may be any interior point on the 
dyadic interval, and can also be an endpoint provided the dyadic
expansion ending in $0$'s is taken at the left endpoint
$\frac{k}{2^n}$ and that ending in $1's$ is taken
at the right endpoint $\frac{k+1}{2^n}.$  

 (2) The functions $\tau_n(x)$ approximate the
Takagi  function monotonically from below
\beql{253b}
\tau_1(x) \le \tau_2(x) \le \tau_3(x) \le ...
\eeq
The values $\{ \tau_n(x):  n \ge 1\}$
converge uniformly to $\tau(x)$, with
\beql{253a}
|\tau_n(x) - \tau(x)| \le \frac{2}{3}\cdot\frac{1}{2^{n}}.
\eeq

(3) For a dyadic rational $x= \frac{k}{2^n}$, perfect approximation
occurs at the $n$-th step, and  
\beql{253c}
\tau(x) = \tau_m (x), ~~~\mbox{for all}~~ m \ge n.
\eeq
\end{theorem}

\paragraph{Proof.} 
All statements follow easily from the observation that 
each function
$f_n(x) := \frac{ \ll 2^n x\gg}{2^n}$
is a piecewise linear sawtooth function, linear on dyadic intervals
$[\frac{k}{2^{n+1}}, \frac{k+1}{2^{n+1}}]$, with slope 
having value $+1$ if the binary expansion of $x$ has $b_{n+1}=0$
and slope having value $-1$ if $b_{n+1}=1$.  The inequality in
\eqn{253a} also uses the fact that $\max_{x \in [0,1]} \tau(x) =
\frac 2 3.$ 
$~~~\Box$\\

%
%

The Takagi function itself can be directly
expressed in terms of the deficient digit function.
The relation \eqn{222} compared with the definition \eqn{203} of 
$\ell_m(x)$ yields
$$
\ell_{m+1}(x) = \frac{m}{2}-  \frac{1}{2} (-1)^{b_{m+1}}\dN_{m}(x).
$$
Substituting this in Takagi's formula \eqn{201} 
and simplifying (noting that    $\ell_{1}(x)= D_0(x)=0$) yields the formula
\beql{223a}
\tau(x)
= \frac{1}{2} -\frac{1}{4} \left(\sum_{m=1}^{\infty} (-1)^{b_{m+1}} \frac{\dN_{m}(x)}{2^m}
\right).
\eeq

%
%
%
%

\section{Basic Properties: Functional Equations and Self-Affine Rescalings} \label{sec4}

We next recall two basic functional equations that the Takagi function satisfies.
These have been repeatedly found, see  Kairies, Darslow and Frank  \cite{KDF88}
and Kairies \cite{Ka98}.

%
%
%
%

\begin{theorem}~\label{th22}
(1) The Takagi function satisfies  two functional equations,
each valid for \\
$0 \le x \le 1,$ the reflection equation
\beql{206a}
\tau(x) = \tau(1-x),
\eeq
and the dyadic self-similarity equation
\beql{206b}
 \tau(\frac{x}{2})=\frac{1}{2} x + \frac{1}{2} \tau(x).
\eeq

(2) The Takagi function on $[0,1]$ is the unique continuous function
on $[0,1]$ that satisfies  both these functional equations.
\end{theorem}

\paragraph{Proof.} 
(1) Here \eqn{206a} follows directly from   \eqn{eq101}, since
$\ll k x\gg ~= ~\ll k(1-x) \gg$ for $k \in \ZZ$.
To obtain \eqn{206b}, let $x= 0. b_1 b_2 b_3 ...$  and set
$y := \frac{x}{2} = 0.0 b_1 b_2b_3 ...$.  Then $\ll y\gg = y$, whence \eqn{eq101} gives
$$
2\tau(y)  = 2 \ll y \gg + 2\left(\sum_{n=1}^{\infty} \frac{ \ll 2^{n} y\gg} {2^n}\right) 
=x + \sum_{m=0}^{\infty}  \frac{ \ll 2^{m} x\gg} {2^m} = x + \tau(x).
$$

(2) The uniqueness result follows by showing that the  functional equations
determine the value at all dyadic rationals. Indeed, the  dyadic self-similarity
equation first gives $\tau(0)=0$, whence $\tau(1)=0$ by the reflection equation.
Then $\tau(\frac{1}{2}) = \frac{1}{2}$ by the self-similarity equation,
and now we iterate to get all other dyadic
rationals. Since the dyadic rationals are dense, there is at most one
continuous interpolation of the function. The fact that a continuous
interpolation exists follows from Theorem \ref{le23}.
(This result was noted by Knuth \cite[Exercise 82, solution p. 103]{Kn05}.
$~~~\Box$\smallskip

The functional equations yield   a self-affine property of the Takagi
function associated to  shifts by dyadic rationals $x = \frac{k}{2^{n}}$.

%
%

\begin{theorem}~\label{le25}
For an arbitrary dyadic rational $x_0=\frac{k}{2^n} $ then for  $x \in [\frac{k}{2^n}, \frac{k+1}{2^n}]$
given by $x = x_0 + \frac{\w}{2^{n}}$, there holds
\beql{251}
\tau( x_0 + \frac{\w}{2^n}) = \tau(x_0) + \frac{1}{2^n} \large( \tau(\w) + \dN_n(x_0)\, \w\large), 
~~~0 \le \w \le 1, 
\eeq
That is, the graph of $\tau(x)$ on  $ [\frac{k}{2^{n}}, \frac{k+1}{2^{n}}]$
is a miniature version of the tilted Takagi function
$ \tau(x) + D_n(x_0) x$, shrunk by a factor
$\frac{1}{2^{n}}$ and
vertically shifted by $\tau(x_0)$. 
\end{theorem}

\paragraph{Proof.} 
By Theorem \ref{le23}(1), we have $\tau_{n}(x_0 + \frac{\w}{2^{n}}) =
\tau_{n}(x_0) + \dN_n(x_0)\cdot \frac{\w}{2^n}$. 
Therefore, by \eqn{eq101} it follows that
\begin{eqnarray*}
  \tau(x) & = & \tau_{n}(x) + \sum_{j=n}^\infty
  \frac{\ll 2^j x\gg}{2^k}\\
  & =  &\tau_{n}(x_0) + \dN_n(x_0)\cdot \frac{\w}{2^n} +  \sum_{j=n}^\infty
  \frac{\ll 2^j (\frac{\w}{2^{n}}) \gg}{2^j}\\
  & =  & \tau(x_0) + \frac{1}{2^n} \large( \tau(\w) + \dN_n(x_0) \w \large).~~~\Box\\
\end{eqnarray*}

Theorem  \ref{le25} simplifies 
 in the special case of $x_0= \frac{k}{2^n}$ having $D_n(x_0)=0$, which 
  we call a  {\em  balanced dyadic rational};
 such dyadic rationals can only occur when  $n=2m$ is  even.  The formula \eqn{251}
 then  becomes 
\beql{251a}
\tau( x_0+ \frac{\w}{2^n}) = \tau(x_0) + \frac{\tau(\w)}{2^n},~~~0 \le \w \le 1,
\eeq
which a  shrinking of the Takagi function together with a vertical shift.
Balanced dyadic rationals play a special role in the  analysis
of the Takagi function.\smallskip

As a final topic in this section, a  convex function is characterized as one that
satisfies the condition
$$
f(\frac{x+y}{2}) \le \frac{f(x) + f(y)}{2}.
$$
The Takagi function is certainly very far from being a convex function.
However one  can establish the following
 approximate mid-convexity property of the
Takagi function, due to Boros \cite{Bor08}

%
%

\begin{theorem}~\label{th43a}
{\rm (Boros 2008)}
The Takagi function is $(0, \frac{1}{2})$-midconvex.
That is, it  satisfies the bound
$$
\tau( \frac{x +y}{2}) \le \frac{ \tau(x) + \tau(y)}{2} + \frac{1}{2} |x -y|.
$$
\end{theorem}

This result establishes the sharpness of general bounds for
approximately midconvex functions established by H\'{a}zy and P\'{a}les \cite{HP04}.

%

\section{The Takagi function and Iteration of  the Tent Map} \label{sec5}

An alternate interpretation of the Takagi function involves 
iterations of  the {\em symmetric tent map}
$T: [0,1] \to [0,1]$, given by 
\begin{equation}\label{eq101a}
T(x) = \left\{ 
\begin{array}{cc} 
2x & \mbox{if}~~ 0\le x \le \frac{1}{2},\\
2-2x & \mbox{if}~~\frac{1}{2}\le x \le 1.
\end{array}
\right.
\end{equation}
This function extends to the whole real line by making it periodic with period $1$,
and is then given by $T(x) = 2 \ll x \gg$,   
and it  satisfies the functional equation $T(x)= T(1-x)$.
A key property concerns its behavior under 
iteration $T^{\circ n}(x) := T( T^{\circ (n-1)}(x))$
($n$-fold composition of functions is here denoted $\circ \, n\,$).

 it satisfies the functional equation under composition
\begin{equation}\label{502a}
T^{\circ n}(x) = T(2^n x).
\end{equation}

The tent map $T(x)$  under iteration is an extremely special map. It
defines under iteration a completely chaotic dynamical system,
defined on the unit interval, whose symbolic dynamics
is the full shift on two letters. This map has Lebesgue measure on $[0, 1]$
as an invariant measure, and this measure is the
maximal entropy measure among all invariant Borel measures for $T(x)$. 

As an immediate consequence of (\ref{502a}) we have the following formula for $\tau(x)$, noted by
Hata and Yamaguti \cite{HY84} in 1984.

%
\begin{theorem}~\label{th44} 
{\rm (Hata and Yamaguti 1984)}
The Takagi function is given by 
\begin{equation}\label{eq101b}
\tau(x) := \sum_{n=1}^{\infty} \frac{1}{2^n} \, T^{\circ n}(x),
\end{equation}
where $T^{\circ n}(x)$ denotes the $n$-th iterate of the tent map $T(x)$.
\end{theorem}

This result   differs conceptually from the original definition \eqn{eq101}
which produces powers of $2$ by a rescaling the variable in a fixed map, in that the powers of $2$ 
 in (\ref{eq101b}) are produced by iteration of a map.

Hata  and Yamaguti \cite{HY84}  defined generalizations of the Takagi
function based on iteration of maps.
They defined the {\em Takagi class} $E_T$  to consist of those  functions
given by 
$$
f(x) = \sum_{n=1}^{\infty} a_n \,T^{\circ n}(x).
$$
where $\sum_{n=1}^{\infty} |a_n| < \infty$
viewed as members of
the  Banach space of continuous functions $C^{0}([0,1])$, under the sup norm.
This class   contains some piecewise
smooth functions, for example
$$
\sum_{n=1}^{\infty} \frac{1}{4^n} T^{\circ n}(x)  = x(1-x).
$$

Hata and Yamaguti showed that continuous functions in $E_T$ can be
characterized in terms of their Faber-Schauder expansions in 
 $C^{0}([0,1])$,
and that the Takagi class $E_T$ is a closed subspace of $C^{0}([0,1])$.
The standard Faber-Schauder  functions, defined by Faber \cite{Fab10} in 1910
and generalized by Schauder \cite[pp. 48-49]{Sch27} in 1927,  are $1, x$ and 
$\{S_{i, 2^j}(x) : j \ge 0, 0\le i \le 2^j-1\}$
 given by  dyadically shrunken and shifted tent maps
$$
S_{i, 2^j}(x) := 2^j \{ |x- \frac{i}{2^{j}}| + |x-\frac{i+1}{2^{j}}| - |2x-\frac{2i+1}{2^{j+1}}|\},
$$
 These functions  form a Schauder basis
  of  the Banach space $C^{0}([0, 1])$,
in the sup norm, 
taking them in the order $1, x$, followed
by  the other functions  in  the order:  $S_{i, 2^j}$ precedes
$S_{i', 2^{j'}}$ if $j < j'$ or $j = j'$ and $i < i'$.
 (For a discussion of Schauder bases see
 Megginson \cite[Sec. 4.1]{Meg98}.)
 A function $f(x)$ in $C^{0}([0,1])$ has 
Faber-Schauder expansion
\begin{equation}\label{eq521}
f(x) = a_0 + a_1 x + \sum_{j=0}^{\infty} \sum_{i=0}^{2^j -1} a_{i, j} S_{i, 2^j}(x),
\end{equation}
with coefficients $a_0= f(0),\,$ $a(1)= f(1)- f(0)$, and 
$$
a_{i, j} := f( \frac{2i+1}{2^{j+1}}) - 
\frac{1}{2}\Big( f(\frac{i}{2^j}) + f(\frac{i+1}{2^j})\Big).
$$
Functions in the Takagi class satisfy the restriction that $a_0=a_1=0$ and their
Faber-Schauder coefficients $\{ a_{i, j}\}$ depend only on the level $j$.

Hata and Yamaguti \cite[Theorem 3.3]{HY84} proved the following
converse result.

%
\begin{theorem}~\label{th52a} 
{\rm (Hata and Yamaguti 1984)}
A function $f(x) \in C^{0}([0,1])$ belongs to the Takagi class $E_T$ if
and only if its Faber-Schauder expansion
$
f(x) = a_0 + a_1 x  +\sum_{j=0}^{\infty} \sum_{i=0}^{2^j-1} a_{i, j} S_{i, 2^j}(x)
$
satisfies 
\begin{enumerate}
\item
The coefficients
$a_0= a_1=0$ 
and 
$$
a_{0, j} = a_{i, j} ~~\mbox{for all} ~~ j \ge 0, \,\, 0  \le i \le 2^{j-1}.
$$
\item If we set $c_j= a_{0, j}$, then 
$$\sum_{j=0}^{\infty} |c_j|< \infty.$$
\end{enumerate}
\end{theorem}

Note that if  the  conditions 1, 2  above hold, then $f(x)$ is given by the expansion
$$
f(x) = \sum_{j=0}^{\infty} c_j \, T^{\circ(j+1)}(x).
$$

Hata and Yamaguti also viewed functions 
in the Takagi  class as  satisfying a 
difference  analogue
of Laplace's equation, using  the scaled   central second difference operators
$$
\Delta_{i, 2^j}(f) := f(\frac{i}{2^j}) + f(\frac{i+1}{2^j})- 2 f( \frac{2i+1}{2^{j+1}}),
$$
for $0\le  i \le 2^{j}-1$, $j \ge 0$, along with Dirichlet
boundary conditions. (The Faber-Schauder coefficients
$a_{i, j} = -\frac{1}{2} \Delta_{i, 2^j}(f)$.)
They obtained the following existence and uniqueness result
(\cite[Theorem 4.1]{HY84}).

%
\begin{theorem}~\label{th53a} 
{\rm (Hata and Yamaguti 1984)}
Given data $\{ c_j: j \ge 0\}$,
the infinite system of linear equations defined for 
$$
\Delta_{i, 2^j}(f) = c_j, ~~~~j \ge 0, ~0 \le i \le 2^j -1, 
$$
has a continuous solution $f(x) \in C^{0}([0,1])$ satisfying Dirichlet
boundary conditions $f(0)=f(1) =0 $ if and only if
$$
\sum_{j=0}^{\infty} |c_j| < \infty.
$$
In this case $f(x) \in E_T$, with
$$
f(x) = -\frac{1}{2} \sum_{j=0}^{\infty} c_j \,T^{\circ(j+1)}(x).
$$
\end{theorem}

There are interesting  functions obtainable from the Takagi
function by monotone changes of variable. 
The  tent map $T(x)$  is  real-analytically conjugate on the interval $[0,1]$ to 
a particular {\em logistic map} 
\begin{equation}\label{eq541a}
F(y) := 4y(1-y), ~~~~~0 \le y \le 1. 
\end{equation}
That is, $F(y)  = \varphi^{-1} \circ T \circ \varphi(y)$ for a monotone increasing
real-analytic function $\varphi(x)$, which is
$$\varphi(y) = \frac{2}{\pi} \arcsin \sqrt{y},$$ 
with its functional inverse $\varphi^{-1}(x)$  given by
$\varphi(x) := \sin^2 \,(\frac{\pi x}{2} ).$
The  dynamics under iteration of the logistic map (\ref{eq541a}) has been much studied.
 It is  a  post-critically finite quadratic polynomial, and 
  it   is affinely  conjugate\footnote{ The conjugacy
is  $y= \psi(z)= -\frac{1}{4} z+\frac{1}{2}$ with inverse $z=\psi^{-1}(y) = -4y+2$.}
to  the monic centered quadratic map $\tilde{F}(z) = z^2 -2,$
which specifies a boundary point of the Mandelbrot set. 

 Its $n$-th  iterate 
$F^{\circ n}(y)$ is a  polynomial of degree $2^n$, which is conjugate
to the Chebyshev polynomial $T_{2^n}(x)$.
The  analytic conjugacy above gives 
  $$
  \varphi \circ F^{\circ n}(y) = T^{\circ n} \circ \varphi(y),  ~~~~\mbox{for}~~n \ge 1.
  $$
 The  change of variable $x= \varphi(y)$ applied
 to a function $f(x) = \sum_{n=1}^{\infty} c_n T^{\circ n}(x)$ 
 in the Takagi class yields  the rescaled function
 \begin{equation}\label{530a}
 g(y) = f(\varphi(y)) = \sum_{n=1}^{\infty} c_n \, \varphi \circ F^{\circ n}(y),
 \end{equation}
motivating the further study of maps of this form.

 More generally, given a dynamical system obtained by
iterating a map $\tilde{T}: [0,1] \to [0,1]$ of the interval, 
and a rescaling function $\psi: [0,1] \to \CC$, 
Hata and Yamaguti define
the {\em generating function}
$$
F(t, x): = \sum_{n=0}^{\infty} t^n \, \psi \circ \tilde{T}^{\circ n}(x).
$$
Here these functions, where $t$ may vary, encode various 
statistical  information about
the discrete dynamical system $\tilde{T}$.
For more information on this viewpoint, see Yamaguti, Hata and Kigami \cite[Chap. 3]{YHK97}.

%
%
%
%

\section{Fourier Series of the Takagi Function}\label{sec6}

The Takagi function defined by \eqn{eq101} extends to 
a continuous periodic function on the real line with period $1$,
which comes with  a Fourier series expansion.
The functional equation $\tau(x)=\tau(1-x)$ for $0 \le x \le 1$ implies that
the extended function also satisfies
\beql{101j}
\tau(x) = \tau(-x),
\eeq
so that it is an even function. As a consequence, its 
Fourier series only involves $\cos (2\pi n x)$ terms. 
It is easily computable from the known  Fourier series
of the symmetric tent map, as follows.

%
\begin{theorem}~\label{th12}
The Fourier series of the Takagi function is given by 
\begin{equation}\label{eq102j}
\tau(x) := \sum_{n \in \ZZ} c_n \,e^{2 \pi i nx},
\end{equation}
in which 
\beql{101k}
c_0  =  \int_{0}^1  \tau(x) dx = \frac{1}{2}
\eeq
and for $n >0$ there holds
\beql{101m}
c_n =  c_{-n} = -\frac{1}{2^{m} (2k+1)^2 \pi^2},  \quad\mbox{where}\quad n= 2^m (2k+1).
\eeq
\end{theorem}

\paragraph{Proof.} The even function $\psi(x)= \ll x \gg$ has Fourier 
series with real coefficients $b_n= b_{-n}$, given as
$$
\psi(x) =  \sum_{n = -\infty}^{\infty} b_n e^{2 \pi i n x} ,
$$
with $b_{-n}=   \int_{-\frac{1}{2}}^{\frac{1}{2}}  |x| e^{2 \pi i n x} dx $.
Clearly $b_0= \frac{1}{4}$ and,
integrating by parts,
\begin{eqnarray*}
\int_{0}^{\frac{1}{2}} x\, e^{2 \pi i n x} dx &=& [  \frac{x}{2 \pi i n} e^{2 \pi in x}] \mid_{x=0}^{x=\frac{1}{2}}
- \int_{0}^{\frac{1}{2}} \frac{1}{2\pi i n} e^{2 \pi i nx} dx\\
&=&  -\frac{1}{4 \pi i n} -[ \frac{1}{ (2 \pi i n)^2} e^{2 \pi i n x} ]\mid_{x=0}^{x=\frac{1}{2}}\\
 &=& \{ \begin{array} {lc} 
 -\frac{1}{4 \pi i n} - \frac{1}{2 \pi^2 n^2}, &~  n ~~\mbox{~odd},\\
 -\frac{1}{4 \pi i n}, & ~ n ~~\mbox{even}.
 \end{array} 
 \end{eqnarray*}
A similar calculation on the interval $[-\frac{1}{2}, 0]$
gives the complex conjugate value, whence 
$$
b_n= b_{-n}=\Big\{\begin{array}{cl}
- \frac{1}{\pi^2 n^2} & ~~\mbox{if  {\em n}~ is odd,}\\
~&~\\
0 &~\mbox{ if}~n \ne 0 ~\mbox{is even}.
\end{array}
$$
Now $\tau(x) = \sum_{m=0}^{\infty} \frac{1}{2^m} \ll 2^m x\gg$, and
the Fourier coefficients of $\psi(2^m x)$ are $b_n^{m} := b_{n/2^m}$
if $2^m | n$ and $0$ otherwise. 
By uniform convergence of the sum 
 we obtain
$c_0 = \sum_{n=0}^{\infty} \frac{1}{2^n} b_0 = \frac{1}{2},$
and, for $n=\pm 2^m(2k+1)$, we obtain
$$c_n = \sum_{j=0}^{\infty}\frac{1}{2^j} b_n^j
=\sum_{j=0}^{m} \frac{1}{2^j} b_{2^{m-j}( 2k+1)} = - \frac{1}{2^{m} (2k+1)^2 \pi^2},$$
which is the result.
$~~~\Box$\\

Note that the decay of the Fourier coefficients
as $n \to \infty$ has $\limsup_{n \to \infty}  n^2 |c_n| > 0.$ 
This fact directly  implies that $\tau(x)$ cannot be a $C^{2}$-function.
However  much more about its oscillator behavior, including
its non-differentiability, can be proved by other methods.

As a direct consequence of this result, the Takagi function is
obtainable as the real part of the boundary value of a holomorphic
function on the unit disk.

%
\begin{theorem}~\label{th32a} 
Let $\{ c_n: n \in \ZZ\}$ be the Fourier coefficients of
the Takagi function, and define the power series
\beql{301a}
f(z) =  \frac{1}{2} c_0 + \sum_{n=1}^{\infty} c_n z^n.
\eeq
This power series converges absolutely 
in the closed unit disk $\{z: |z| \le 1\}$ to define a holomorphic function
in its interior,  and it has the unit circle as a natural boundary
to analytic continuation.
It defines a continuous function on the boundary of the unit disk, 
and its values there satisfy
\beql{302a}
f(e^{2 \pi i \theta}) = \frac{1}{2} \left(X(\theta) + i Y(\theta)\right)
\eeq
in which $X(\theta)=\tau(\theta)$ is the Takagi function. \end{theorem}

\paragraph{Remark.} The imaginary part defines a
function $Y(\theta)$ is a new function which
we term the {\em conjugate Takagi function.}

\paragraph{Proof.} The Fourier coefficients 
of the Takagi function satisfy
\[
\sum_{n=0}^{\infty} |c_n|= \frac{1}{2} +
\sum_{m=0}^{\infty} \frac{1}{{2^{m-1}}}\left( \sum_{k=0}^{\infty} \frac{1}{(2k+1)^2\pi^2}\right)
=\frac{1}{2} + \frac{1}{2} < \infty
\]
It follows that the power series for $f(z)$ converges absolutely on the unit circle,
so is continuous on the closed unit disk, and holomorphic in its interior.
Since the Fourier series for the Takagi function is even, by inspection
\[
f(e^{2 \pi i \theta}) + f(e^{- 2 \pi i \theta}) = 2 Re( F(e^{2 \pi i \theta})) = X(\theta).
\]
 This justifies \eqn{302a}, defining $U(\theta)$ by
\[
f(e^{2 \pi i \theta}) - f(e^{- 2 \pi i \theta}) = 2 i Im( F(e^{2 \pi i \theta})) = i Y(\theta).
\]
Since the Takagi function is non-differentiable everywhere, the function $f(z)$
cannot analytically continue across any arc of the unit circle, so that the 
unit circle is a natural boundary for $f(z)$.
$~~~\Box$\\

As mentioned in Section \ref{sec5}, the  Takagi function can also be
studied using its Faber-Schauder expansion, rather than a Fourier expansion.
In 1988 Yamaguti and Kigami  \cite{YK88}  defined for $p \ge 1$ the Banach spaces
$$
H_p := \overline{ \{ f(x) = \sum_{i,j} c_{i, j}  2^{-\frac{j}{2} +1}S_{i, 2^j}(x) \} },
$$
defined using a rescaled Schauder basis and the  norm $||f||_p = \Big( \sum_{i,j} |c_{i,j}|^p \Big)^{1/p}.$
They deduced using its Schauder
expansion  that the Takagi  function belongs to the Banach space $H_p$ for all $p>2$
(It does not belong to $H_2$, which is a Hilbert space coinciding with a space denoted $H^{1}$ in \cite{YK88}.)
 In 1989 Yamaguti \cite{Yam89} proposed a generalized  Schauder basis consisting
 of polynomials in $x$, 
in which to study the Takagi function and other functions.

%
%
%
%

\section{The Takagi Function and
Bernoulli Convolutions} \label{sec7}

The Takagi function   also appears in
the analysis of non-symmetric Bernoulli convolutions. Let
$$
x = \sum_{j=1}^{\infty} \epsilon_j (\frac{1}{2})^j,
$$
in which the digits  $\epsilon_j \in \{0,1\}$ are drawn as independent  (non-symmetric)
Bernoulli random variables taking
 value $0$ with probability $\alpha$ and $1$ with probability $1- \alpha$.
Then $x$ is a random real number in $[0,1]$ with cumulative distribution function
$$
L_{\alpha}(x) = \mu_{\alpha}([0,x])= \int_{0}^x d\mu_{\alpha},
$$
in which $\mu_\alpha$ is a  certain  Borel measure on $[0,1]$.  
These functions were introduced in 1934 by Lomnicki and Ulam \cite[pp. 267-269]{LU34}
 with this interpretation.
The measure $\mu_{\alpha}$  is Lebesgue measure for
$\alpha= \frac{1}{2}$, and is a singular measure otherwise. \\

In 1943  R. Salem \cite{Salem43} gave a geometric construction of
monotonic increasing singular functions  that includes these functions as a special case. 

%
\begin{theorem}{\rm (Salem 1943)}\label{th21} 
The function $L_{\alpha}(x)$ has  Fourier series 
 $L_\alpha(x) \sim \sum_{n=1}^{\infty} c_n \cos 2 \pi nx$
using the formula (for real $t$) 
$$
\int_{0}^{1} e^{2\pi itx} dL_{\alpha}(x) = \prod_{k=1}^{\infty}\left( \alpha + (1-\alpha)
e^{\frac{2\pi i t}{2^k}}\right).
$$
The  Fourier coefficients are given by infinite products 
$$
c_n = e^{- \pi i n} \prod_{k=1}^{\infty} \left(  \cos \frac{\pi n}{2^k} + i(2 \alpha -1)\sin \frac{\pi n}{2^k} \right).
$$
\end{theorem}

An interesting property of the
Lomnicki-Ulam  function $L_{\alpha}(x)$  
is that it satisfies a {\em two-scale dilation equation}
$$
L_{\alpha}(x)=\left\{ 
\begin{array}{cl}  \alpha L_{\alpha}(2x) & ~\mbox{for} ~~0 \le x \le \frac{1}{2},\\
 (1-\alpha)L_{\alpha}(2x-1) + \alpha & ~\mbox{for} ~~\frac{1}{2} \le x \le 1.
 \end{array} 
 \right.
 $$
In 1956 de Rham (\cite{deR56} , \cite{deR56b}, \cite{deR57}) 
studied such functional equations in detail.
Dilation equations are relevant to the 
construction  of compactly
supported wavelets,  and general properties of solutions to such
equations were derived  in 
Daubechies and Lagarias \cite{DL92a}, \cite{DL92b}.

In 1984 Hata and Yamaguti \cite[Theorem 4.6]{HY84} made the following connection of 
the Lomnicki-Ulam functions $L_{\alpha}(x)$ to the Takagi function, cf. \cite[p. 195]{HY84}. 
%
\begin{theorem} {\rm (Hata and Yamaguti (1984))} \label{th22a}
For fixed $x$ the function
$g_x(\alpha) := L_{\alpha}(x)$, initially defined for $\alpha \in  [0,1]$, extends to an analytic function 
of $\alpha$ on the
lens-shaped region 
$$
D=\{ \alpha \in \CC:~ |\alpha| < 1 ~~\mbox{and}~~|1- \alpha| < 1\}.
$$
The Takagi function
appears as the derivative of these functions at the central point $\alpha= \frac{1}{2}$: 
\beql{111}
\frac{d}{d\alpha} L_{\alpha}(x)|_{\alpha= \frac{1}{2}} =  2 \,\tau(x).
\eeq
\end{theorem} 

This result permits an interpretation of the Takagi function as a generating function
of a chaotic dynamical system (Yamaguti et al. \cite[Chapter 3]{YHK97}).
This result was generalized further by Sekiguchi and Shiota \cite{SS91}.

%
%
%
%

\section{Oscillatory Behavior  of the Takagi Function} \label{sec8}

The main feature of the Takagi function is that it is non-differentiable
everywhere. Takagi \cite{Tak03} showed the function has no two-sided finite
derivatives at any point.  In 1984 Cater \cite{Ca84}
showed that the Takagi function has no one-sided finite derivative at any point.
However it does have well-defined (two-sided) improper derivatives equal to $+\infty$
(resp. $-\infty$) at some points, see Theorem \ref{th86a}. 

The non-differentiability of the Takagi function is bound up with its
increasing oscillatory behavior as the scale decreases.   
Letting $0 < h <1$ measure
a scale size,  we have
the following elementary estimate bounding the maximal size of
oscillations at scale $h$.

%
\begin{theorem}\label{th91a}
For $0\le x \le x+h \le 1$, the Takagi function satisfies
\beql{601a}
|\tau(x+h) - \tau(x)| \le 2 h \log_2 \frac{1}{h}.
\eeq
\end{theorem}

\paragraph{Proof.}
Suppose $2^{-n} \le h \le 2^{-n+1}$, so that $n   \le \log_2 \frac{1}{h}.$
Theorem ~\ref{le23} gives the estimate
that $|\tau(x) - \tau_n(x)| \le \frac{2}{3}\frac{1}{2^n}$,
 and  we know $\tau_n(x)$ has 
 everywhere slope between $-n$ and $n$.  It follows that
 $$
 |\tau(x+h) - \tau(x)| \le \frac{2}{3}\left( \frac{1}{2^n}\right) + n h \le h(n+\frac{2}{3}) \le 2h \log_2\frac{1}{h}.
 $$
as required.$~~~\Box$.\medskip

Theorem \ref{th91a} is sharp to within a multiplicative factor of $2$, since for $h=2^{-n}$,
$$
\tau(h) - \tau(0) = \tau(2^{-n}) = \frac{n}{2^n} = h \, \log_2 \frac{1}{h},
$$
In fact the  multiplicative factor of $2$ can be decreased  to $1$ as
the scale decreases,  in the following sense (K\^{o}no \cite[Theorem 4]{Ko87}).

%
\begin{theorem}\label{th92a} {\rm (K\^{o}no 1987)}
Let $\sigma_u(h) =  \log_2 \frac{1}{h}$. Then there holds
\[
\limsup_{|x-y| \to 0^{+}} \frac{\tau(x)- \tau(y)}{|x-y| \, \sigma_u(|x-y|)} =1.
\]
and
\[
\liminf_{|x-y| \to 0^{+}} \frac{\tau(x)- \tau(y)}{|x-y| \, \sigma_u(|x-y|)} = -1.
\]
\end{theorem}

The average size of extreme fluctuations for most $x$  is of the smaller order 
$h \sqrt{\log_2 \frac{1}{h}} \sqrt{2 \log\log \log_2 \frac{1}{h}}$,
as given in the following result.
(K\^{o}no \cite[Theorem 5]{Ko87}).

%
\begin{theorem}\label{th93a} {\rm (K\^{o}no 1987)}
Let $\sigma_l(h) = \sqrt{ \log_2 \frac{1}{h}}$.
Then for almost all $x \in [0,1]$ there holds
\[
\limsup_{h \to 0^{+}} \frac{\tau(x+h)- \tau(x)}{h \,\sigma_l(h) \sqrt{2 \log\log \sigma_l(h)}} =1,
\]
and
\[
\liminf_{h \to 0^{+}} \frac{\tau(x+h)- \tau(x)}{h \,\sigma_l(h) \sqrt{2 \log\log \sigma_l( h)}} = -1.
\]
\end{theorem}

K\^{o}no used expansions of $\tau(x)$ in terms of Rademacher functions
to obtain his results.\smallskip

Finally, if one scales the oscillations by the factor
$h \sqrt{ \log_2\frac{1}{h}}$ on the scale size $h$ then one obtains a Gaussian limit
distribution of the maximal oscillation sizes at scale $h$ as $h \to 0^{+}$, in the following sense
(Gamkrelidze\cite[Theorem 1]{Ga90}).

%
\begin{theorem}\label{th94a}{\rm (Gamkrelidze  1990)}
Let $\sigma_l(h) = \sqrt{ \log_2 \frac{1}{h}}$. Then for each real $y$, 
\beql{611a}
\lim_{h \to 0^{+}} \, Meas \,\{x: \frac{\tau(x+h) - \tau(x)}{h \,\sigma_l(h)} \le y\} =
\frac{1}{\sqrt{2 \pi}} \int_{-\infty}^y e^{-\frac{1}{2} t^2} dt.
\eeq

\end{theorem}

\noindent Gamkrelidze also
observes that Theorem \ref{th93a} can be derived as a consequence of this
result. Namely,  in this context Theorem \ref{th93a}
is interpretable as analogous to the law of the
iterated logarithm in probability theory.\\

The oscillatory behavior of the Takagi function has also been
studied in the H\"{o}lder sense. 
We define on $[0,1]$ the class  $C^0$ of continuous functions, the class 
$C^1$ of continuously differentiable functions (with one-sided
derivatives at the endpoints)  and,
for $0 < \alpha \le 1$, the (intermediate) Lipschitz classes
\beql{511}
Lip^{\alpha} := \{ f\in C^0:~ \mbox{there exists} ~K >0 ~\mbox{with} ~ |f(x) - f(y)| < K|x-y|^{\alpha},
~x,y \in [0,1]\}.
\eeq

%
\begin{theorem} ~\label{th95a} {\rm (de Vito 1985;~ Brown and Kozlowski 2003)}\\
(1) The Takagi  function $\tau$ belongs to the function class
\beql{512}
\tau \in \bigcap_{0< \alpha < 1} Lip^{\alpha}.
\eeq

(2) The Takagi function
 $\tau$ does not agree with any function $g \in C^1$ on any set of positive measure.
In fact, if $M$ is a subset of $[0,1]$ with positive measure then the set
\beql{513}
D(\tau, M) := \{ \frac{\tau(y)-\tau(x)}{y-x}: ~~ x, y \in M ~~\mbox{with}~~x \ne y\}
\eeq
is unbounded. 
\end{theorem}

\paragraph{Proof.} (1) This was shown by de Vito \cite{dV58} in 1958.
He showed moreover, that for $0 < \alpha <1$ and $0 \le x, y \le 1$, there holds
$$
|\tau(x) - \tau(y)| \le \frac{2^{\alpha-1}}{1-2^{\alpha-1}} |x- y|^{\alpha}.
$$
Another proof of  \eqn{512} was given by Shidfar and Sabetfakhri \cite{SS86}.
An extension to $\tau_r$ for all even $r \ge 2$ in (\ref{eq102a}) follows from results in
Shidfar and Sabetfakhri \cite{SS90}. 

(2) A proof  establishing this and  \eqn{513}
was given by Brown and Kozlowski \cite{BK03} in 2003. $~~~\Box$\\

Concerning the non-differentiability of the Takagi function, 
Takagi established there is no finite derivative at any point of $[0,1]$.
However  the Takagi function does have some points where it has a 
well-defined ( two-sided)  infinite deriviative.
These points were recently classified by Allaart and Kawamura \cite{AK10},
who proved the following result (\cite[Corollary 3.9]{AK10}).

%
\begin{theorem} ~\label{th86a} {\rm (Allaart and Kawamura 2010)}
The set of points where the Takagi function has a well-defined
two-sided derivative $\tau^{'}(x) = +\infty$, and the set of points
where it has $\tau^{'}(x) = -\infty$ both are dense in $[0,1]$
and have Hausdorff dimension $1$.
\end{theorem}

There has been much further work studying oscillatory behavior of
the Takagi function in various metrics. 
We mention the work of  Buczolich \cite{Buz03} and  
of Allaart and Kawamura \cite{AK06}, \cite{AK10}.
In addition Allaart \cite{Aal09} has proved
analogues of many results above for a wider class
of functions in the Takagi class.

%
%
%
%

\section{The Takagi function and Binary Digit Sums}\label{sec9}

The Takagi function appears  in the analysis of   binary digit 
sums. 
More generally,
let $A_r(x)$ denote the sum of the base $r$ digits of all integers below $x$, so that 
\beql{112}
A_r(x) := \sum_{j=1}^{[x]} S_r(j),
\eeq
in which $S_r(j)$ denotes the  
sum of the  digits in the base $r$  expansion of $j$. In 1949 Mirsky \cite{Mir49}
observed that 
$$
A_r(x) = \frac{1}{2} (r-1)  x \log_r x - E_r(x)
$$
with remainder term $E_r(x) = O(x)$.
In 1968 Trollope \cite{Trollope68} observed that, for $r=2$,
this  remainder term has an explicit formula given in terms of the
Takagi function. 


\begin{theorem} {\rm (Trollope 1968)}\label{th71a} 
If  the integer $n$
satisfies $2^m \le n < 2^{m+1}$, and if we write $n = 2^m(1+x)$
for a rational number $0 \le x< 1$, with $x= \frac{n}{2^m} -1$,  then
$$
E_2(n) = 2^{m-1}\left( (1+x) \log_2 (1+x) - 2x + \tau(x) \right).
$$
\end{theorem}

This was generalized in a very influential paper of Delange \cite{Delange75},
who observed the following result.
%
\begin{theorem} {\rm (Delange 1975)} \label{th72}
One can write the rescaled error term $E_r(n)$ for base $r$ digit sums at integer values $n$  as 
$$
\frac{1}{n} E_r(n) = \frac{1}{2} \log_r n + F_r \left( \log_r n\right),
$$
in which $F_r(t)$ is a continuous  function  which is periodic of period $1$.
The Fourier series of the function $F_r(t)$ is given by
 $$
 F_r(t)= \sum_{k \in \ZZ} c_k (r) e^{2 \pi i k t}, 
 $$ 
The
Fourier coefficients $c_k(r)$
 involve the values of the Riemann zeta function $\zeta(s)$ at the
points $\frac{2\pi k i}{\log r}$ on the imaginary axis. 
\end{theorem}

In view of Trollope's result,
in  the case $r=2$ the Takagi function $\tau(x)$ appears in the function $F_2(t)$,
namely this function is 
\beql{715}
F_2(t) = \frac{1}{2}( 1+ \lfloor t\rfloor - t) + 2^{1+ \lfloor t\rfloor -t} \tau(2^{t - \lfloor t\rfloor -1}).
\eeq
Also for $r=2$,
and  $k \ne 0$,  the Fourier coefficients of $F_2(t)$ are
$$
c_k= \frac{\frac{i}{2 k \pi}}{1+ \frac{2 k \pi i}{\log 2}}\, \zeta(\frac{2 k \pi i}{\log 2}),
$$
and for $k=0$, 
$$
c_0 = \frac{1}{2\log 2}( \log 2 \pi -1) - \frac{3}{4}.
$$

Results on $k$-th powers  of binary digit sums were
given in 1977 by Stolarsky \cite{St77}, who also provides a summary of earlier literature.
A  connection of these sums with a series of more complicated oscillatory functions 
was obtained in 1986  in Coquet \cite{Coq86}.
Further work was done by Okada, Sekiguchi and Shiota \cite{OSS95} 
in 1995. Recently Kr\"{u}ppel \cite{Kr08} gave another derivation of these
results and further extended them. \smallskip

For a related treatment of similar functions, 
using Mellin transforms, see 
Flajolet et al \cite{FGKPT94}. 
A general survey of dynamical systems in numeration, including
these topics, 
is given by Barat, Berth\'{e}, Liardet and Thuswaldner \cite{BBLT06}.

%
%
%
%

\section{The  Takagi Function and the Riemann Hypothesis} \label{sec10}

The Riemann hypothesis can be formulated purely  in terms of the Takagi function,
as was observed by Kanemitsu and Yoshimoto \cite[Corollary to Theorem 5]{KY00}.
Their  encoding of the Riemann hypothesis
concerns the values of the function restricted to the Farey
fractions $\sF_N$. \smallskip

The Farey series of level $N$ consists of all
reduced rational fractions $0 \le \frac{p}{q} < 1$ having denominator at most $N$,
their number being $|\sF_N| = \frac{3}{\pi^2} N^2 + O( N \log N)$.
The connection of the approximately uniform 
 spacing of Farey fractions and the Riemann hypothesis
starts with  Franel's theorem, cf. Franel \cite{Fr24}, proved in 1924.
(cf. Huxley \cite[p. 36]{Hu72}). Since that time many variants of Franel's result have
been established.  A particular version  of Franel's theorem that is relevant here
was given by Mikolas \cite{Mik49}, \cite{Mik50}, \cite{Mik51}.\smallskip

Uniformity  of distribution of a set of points (here the
Farey fractions)  can be  measured in terms of the efficiency
of numerical integration of functions on the unit interval obtained by sampling
their values at these points. This is the framework taken
in the result of Kanemitsu and Yoshimoto. 
The particular interest of their result
is that it applies to the numerical integration of continuous
functions that are not necessarily  differentiable, but which
belong to the Lipschitz class $Lip^{1- \epsilon}$ for
each $\epsilon >0$. The role of the Takagi function here
is to be a particularly interesting example  to
which their general theorem applies.

%
\begin{theorem} ~\label{th81a}
{\rm (Kanemitsu and Yoshimoto 2000)}
The Riemann hypothesis is equivalent to the statement that for each $\epsilon>0$
there holds, for the Farey sequence $\sF_N$, 
\beql{811b}
\sum_ {\rho \in \sF_N} \tau(\rho) - |\sF_N| \int_{0}^1 \tau(x) dx  = O\left( N^{\frac{1}{2} + \epsilon}\right)
\eeq
as $N \to \infty$. 
\end{theorem}

\paragraph{Proof.}
This result was announced in Kanemitsu and Yoshimoto \cite[Corollary to Theorem 5]{KY00},
and a detailed proof  was given in 2006 in
Balasubramanian, Kanemitsu and Yoshimoto \cite{BKY06}. 
According to the Principle stated in \cite[p. 4]{BKY06},
 the exponent $\frac{1}{2}+\epsilon$ in the remainder term
in \eqn{811b} depends on the fact that the
Takagi function is in the Lipschitz class $Lip^{\alpha}$ for any $\alpha < 1,$
see Theorem ~\ref{th94} below.
(Note that  in \cite[Sec. 2.2]{BKY06} their Takagi function $T(x)$ is defined by a Fourier series,
which when compared with  Theorem ~\ref{th12} suggests it needs
 a rescaling plus a constant term added to agree with $\tau(x)$; this does not
affect the general argument.)
$~~~\Box$\medskip

An interesting feature
of this result is that  for each $N$ the left side of \eqn{811b} is a rational number.
This follows since  the Takagi function takes rational values at
rational numbers, and because with our scaling of  the Takagi function we  have  $\int_{0}^1 \tau(x) dx=\frac{1}{2}$.


\section{ Graph of the Takagi Function} \label{sec11}

We next consider properties of  the graph of the Takagi function
$$
\sG(\tau):= \{ (x, \tau(x)): ~0 \le x \le  1 \}.
$$

The  local extreme points of the graph of the Takagi function were
determined by Kahane \cite{Kah59}.

%
\begin{theorem} ~\label{th83} 
{\rm  (Kahane 1959) }
 
 (1) The  set of local minima
 of the Takagi function are exactly the set of all 
dyadic rational numbers in $[0,1]$.

(2) The set of local maxima of the Takagi function are exactly those
points $x$ such that the binary expansion of $x$ have 
deficient digit function $\dN_{2n}(x)=0$  for all sufficiently large $n$.
\end{theorem}

\paragraph{Proof.} This is shown in \cite[Sec. 1]{Kah59}.
Kahane states condition (2) as requiring that the binary
expansion of $x$ have $b_1+ b_2+ \cdots + b_{2n} = n$ holding
for all sufficiently large $n.$ $~~~\Box$\smallskip

There are uncountably many $x$ that satisfy condition (2) above.
It is  easy to deduce from 
this result that 
the graph of the Takagi function contains a dense set of local minima and local maxima,
viewed from either of the abscissa or ordinate directions.
Allaaert and Kawamura \cite{AK06} further study  the extreme values of functions 
related to the Takagi function.\smallskip

The Hausdorff dimension of the graph of the
Takagi function was determined by Mauldin and Williams \cite{MW86}.

%
\begin{theorem} ~\label{th81} {\rm (Mauldin and Williams 1986)}
 The graph $\sG(\tau)$ of the Takagi function has Hausdorff dimension $1$
as a subset of $\RR^2$. 
\end{theorem}

\paragraph{Proof.} This result is given 
 as \cite[Theorem 7]{MW86}.
 $~~\Box$ \medskip

Mauldin and Williams \cite{MW86} raised the question of
whether the graph of the Takagi function is $\sigma$-finite. 
This was answered in the affirmative by Anderson and Pitt 
\cite{AP89}. 

%
\begin{theorem} ~\label{th82} {\rm (Anderson and Pitt 1989)}
The graph $\sG(\tau)$ of the Takagi function has  $\sigma$-finite linear measure.
\end{theorem}

 \paragraph{Proof.} This result is given 
 as  \cite[Thm. 6.4, and Remark p. 588]{AP89}.
 $~~\Box$ \medskip
 
In the opposite direction, one can show the graph of the Takagi function
has infinite length. We have the following result.

%
\begin{theorem} ~\label{th82}
The graph $\sG(\tau)$ has infinite length locally. That it, it has
infinite length over any
nonempty open  interval $x_1 < x < x_2$ in $[0,1]$.
\end{theorem}

\paragraph{Proof.}This is proved using the piecewise linear 
approximations $\tau_n(x)$ to the Takagi function.
This function has control points at $x=\frac{k}{2^n}$ where it
takes  values agreeing with $\tau(x)$. It suffices to show 
that the length of $\tau_n(x)$ becomes unbounded as $n \to \infty.$ This shows
the whole graph has infinite length. 
The self-similar functional equation then implies that
any little piece of the graph over the interval $x_1 \le x \le x_2$ also  has infinite length.

To show the unboundedness of the length of $\tau_n(x)$, as $n \to \infty$, it suffices to show
that the average size of the (absolute value) of the slope of $\tau_n(x)$ becomes unbounded.
To do this we note that the subdivision from level $n$ to level $n+1$ replaces
each slope $m$ with two intervals of slopes $m +1$ and $m-1$. These do not
change the average value of the (absolute value) of slope, except when $m=0$.
Such slopes occur only on even levels $2n$ and there are $\Big( {{2n}\atop{n}}\Big)$
of them (they are the balanced dyadic rationals at that level). Consequently the
average value of the slope at level $2n$ (resp. $2n+1$) obeys the recursion:
$a_{2n}= a_{2n-1}$ and $a_{2n+1} = a_{2n} + \frac{1}{2^{2n}} \Big( {{2n}\atop{n}}\Big)$.
Since 
$$\frac{1}{2^{2n}}\Big( {{2n}\atop{n}}\Big)\sim \frac{1}{\sqrt{4 \pi n}}$$
and $\sum_{n=1}^{\infty} \frac{1}{\sqrt{4\pi n}}$ diverges, we have
$a_n \to \infty$ as $n \to \infty$, proving the result. $~~~\Box$.\medskip

In 1997 Tricot \cite{Tri97}
defined the notion of {\em irregularity degree} of a function $f$ to be
the Hausdorff dimension of its graph $\sG(f)$. He introduces a two-parameter
family of norms $\Delta^{\alpha, \beta}$ to measure oscillatory behavior.  
As an example, he computes these values for
 the Takagi function and its relatives (\cite[Cor. 6.5, Cor. 6.7]{Tri97})

%
%
%

\section{Level Sets of the Takagi Function} \label{sec12}

We next  consider properties of  level sets of the Takagi function. 
For $y \in [0, \frac{2}{3}]$  we denote the
{\em  level set}  at level $y$   by 
$$
L(y):= \{ x: ~ \tau(x) =y, ~~0 \le x \le 1\}. 
$$
The level  sets have a complicated  and interesting structure.

%
\begin{theorem} ~\label{th61} \label{th61a}
(1) The Takagi function has level sets $L(y)$ that take a finite number of values,
 resp. a countably infinite set of values, resp. an
uncountably infinite set of values.
Specificially, the level set  $L(0)= \{ 0, 1\}$ takes two 
values, the level set  $ L(\frac{1}{2})$ takes a countably infinite
set of values, and the level set $L(\frac{2}{3})$ takes an
uncountably infinite set of values.

(2) The set of levels $y$ such that $L(y)$ is  finite, resp. countably infinite,
resp uncountably infinite,  
are each dense in $[0, \frac{2}{3}]$.
\end{theorem}

\paragraph{Proof.}
(1) It is clear that $L(0)=\{0,1\}$ is a finite level set.
A finite example given by Knuth \cite{Kn05} is $L(\frac{1}{5}) = \{ x, 1-x\}$
with $x= \frac{83581}{87040}.$  Knuth \cite[pp. 20-21, 32--33]{Kn05}
also showed that
certain dyadic rationals $x= \frac{m}{2^k}$ have countably infinite level sets;  these include $x= \frac{1}{2}.$
 The author  and Maddock \cite[Theorem 7.1]{LM10a} also show that $L(\frac{1}{2})$
 is countably infinite. For uncountably infinite level sets, 
Baba \cite{Baba84} showed  that $L(\frac{2}{3})$ has positive Hausdorff dimension,
which implies it is uncountable.

(2) For uncountably infinite sets, this follows from the self-similarity
relation in Theorem \ref{le25}. For countably infinite case it follows from
Knuth's results, see also Allaart \cite{A11b}. For the finite case, it follows
from Theorem ~\ref{th94} below.
$~~~\Box$.\medskip

We next consider bounds on the Hausdorff dimension of level sets. 

%
\begin{theorem} ~ \label{th92} {\rm (Baba 1985; de Amo et al 2011)}

(1) The maximal level set $L(\frac{2}{3})$ has Hausdorff dimension $\frac{1}{2}$.
The set of  levels having level set of
 Hausdorff dimension $\frac{1}{2}$ is dense in $[0, \frac{2}{3}]$.\medskip

(2) Every level set $L(y)$ of the Takagi function has Hausdorff dimension
at most $\frac{1}{2}$. 
\end{theorem}

\paragraph{Proof.} 
(1) In 1984 Baba \cite{Baba84}
 showed  that $L(\frac{2}{3})$ has Hausdorff dimension $\frac{1}{2}$.
Using the self-similarity Theorem ~\ref{le25} for balanced dyadic rationals,
the same holds for $y= \tau(x)$ such that $x$ has a binary expansion
whose purely periodic
part is $(01)^{\infty}$, and whose preperiodic part has an equal number of zeros
and ones. The set of such $x$ is dense in $[0,1]$, so their
image values $y$ are dense in $[0, \frac{2}{3}]$.

(2) An upper bound on the Hausdorff dimension of $0.699$ was
established in 2010 by Maddock \cite{Mad10}. Recently de Amo, Bhouri, D\'{i}az Carrillo and Fern\'{a}ndex-S\'{a}nchez
\cite{ABDF11}
improved on his argument to establish the optimal upper bound of $\frac{1}{2}$
on the Hausdorff dimension.
$~~~\Box$.\smallskip

\begin{figure}[h]
  \begin{center}$                                                               
    \begin{array}{cc}
      \includegraphics[height=2in]{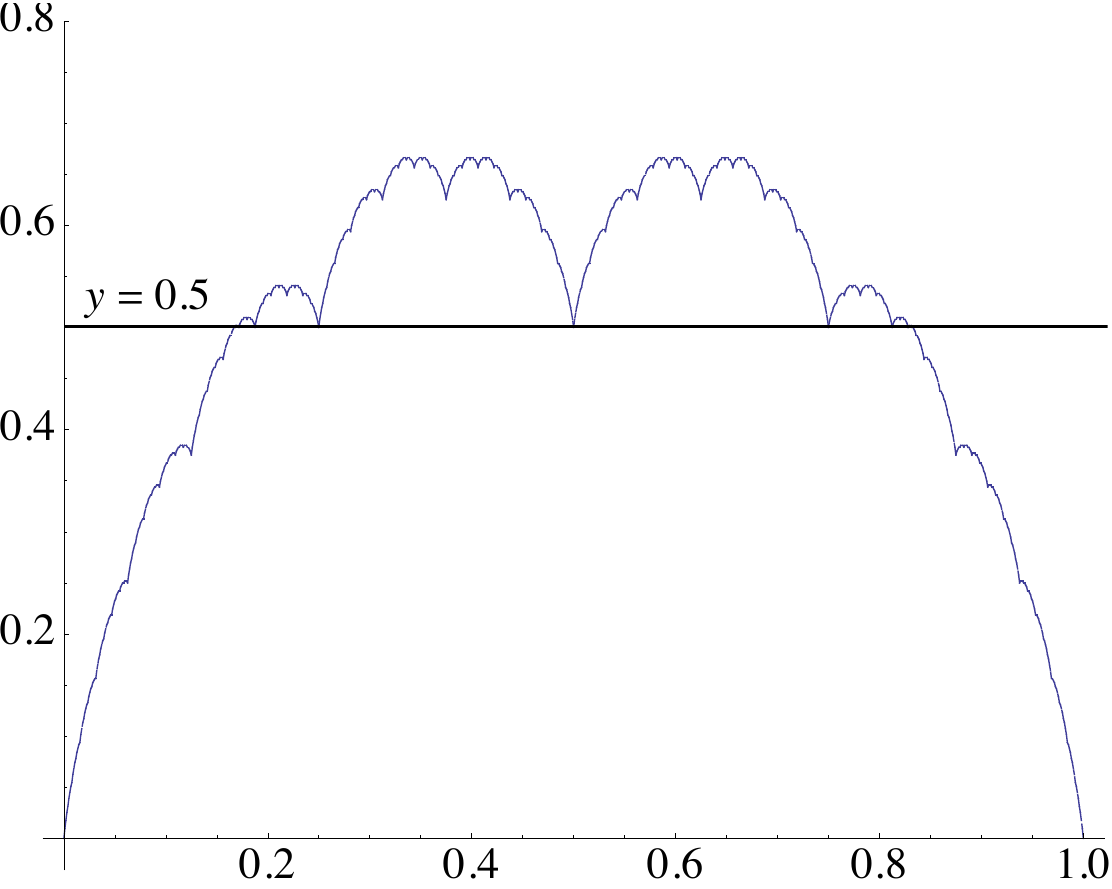} &
      \includegraphics[height=2in]{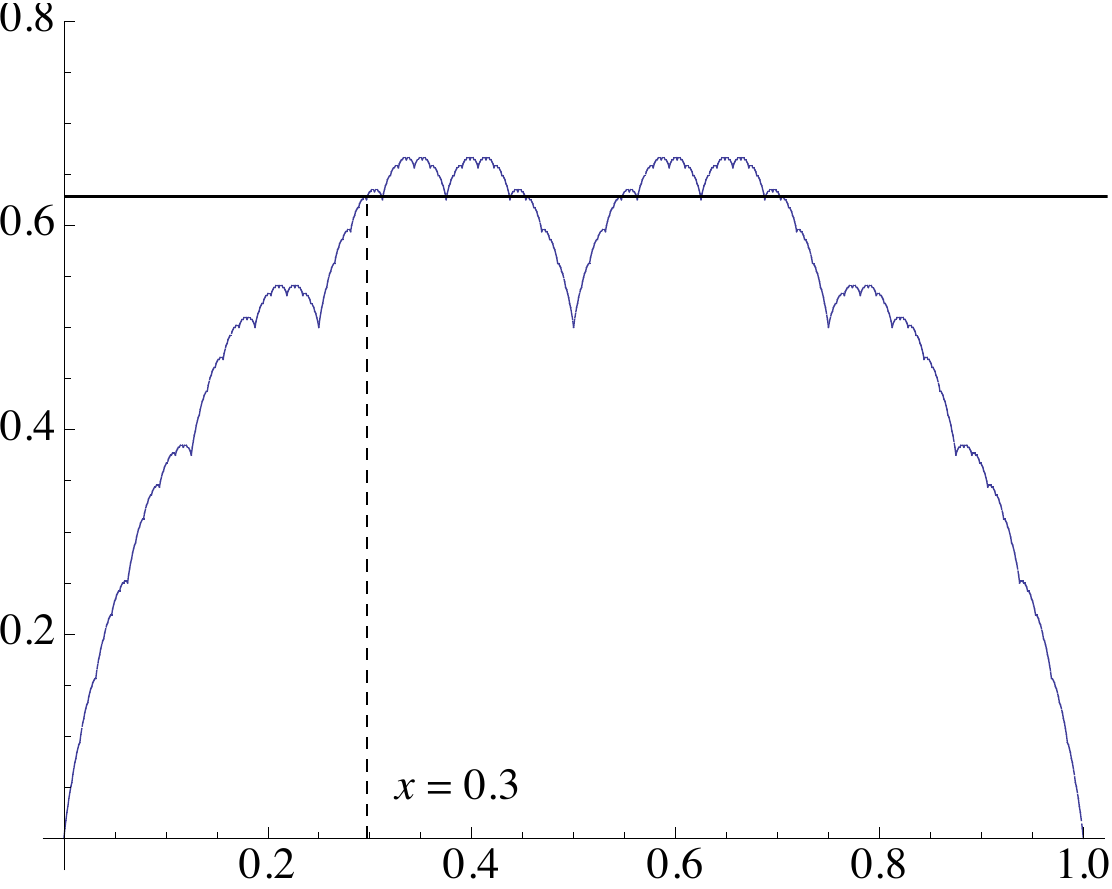}
    \end{array}$
  \end{center}
\caption{Ordinate level  set $L(y)$ at  $y=0.5$ and  abscissa level set
  $L(\tau(x))$ at   $x=0.3$.}
\label{fig111}
\end{figure}

We next consider the nature of  ``generic level sets", considered
in two different senses.  The first is, to draw a abscissa value  $x$ at
random in $[0,1]$ with respect to Lebesgue measure, and to ask about
the nature of the level set $L(\tau(x))$. The second is, to draw an ordinate value 
$y$ at random in $[0,  \frac{2}{3}]$, and ask what is the nature of the level
set $L(y)$. Figure \ref{fig111} illustrates the two senses.\medskip

%
%

For sampling random abscissa values the result is as follows.
%
%
\begin{theorem}\label{th93}
For a full Lebesgue measure set of abscissa points $x \in [0,1]$ the level set
$L(\tau(x))$ contains an uncountable Cantor  set. 
\end{theorem}

\paragraph{Proof.}
This is a corollary of a result of the author and  Maddock \cite[Theorem 1.4]{LM10a},
The latter result concerns local level sets,  discussed below in Section \ref{sec13}.
$~~~\Box$. \smallskip

For sampling random ordinate values the corresponding result is quite different.
%
%
\begin{theorem}~\label{th94} {\rm (Buczolich 2008)}
 For a full Lebesgue measure set of ordinate points $y \in  [0, \frac{2}{3}]$ the
 level set $L(y)$ is a finite set. 
 \end{theorem}
 
 \paragraph{Proof.} This result was proved by Buczolich \cite{Buz08} in 2008.$~~\Box$\medskip
 
 A refinement of the ordinates result is as follows.
%
%
\begin{theorem}~\label{th94a}
Letting $|L(y)|$
denote the number of elements in $L(y)$, if finite, and setting $|L(y)|=0$ otherwise,
one has 
$\int_{0}^{\frac{2}{3}} |L(y)| dy = + \infty.$
That is, the expected number of elements in 
$|L(y)|$
for $y \in  [0, \frac{2}{3}]$  drawn uniformly,  is infinite. 
\end{theorem}

\paragraph{Proof.} 
This result was obtained in \cite[Theorem 6.3]{LM10a}.
A simplified proof was given by Allaart \cite{A11}.
$~~\Box$\medskip

The ordinate and abscissa results 
on the size of level sets are not contradictory; sampling a point $x$ on the
abscissa will tend to pick level sets which are ``large". In order for these two
results to hold 
it appears necessary that the Takagi function must (in some imprecise sense)
have  ``infinite slope" over part of its domain.

The following result goes some way
towards reconciling these two results, 
by showing that the set of large ``ordinate" level sets has
full Hausdorff dimension.

%
%
\begin{theorem}~\label{th95} 
Let $\Gamma_H^{ord}$ be  set of ordinates $y \in [0, \frac{2}{3}]$ such that the
Takagi function  level set $L(y)$ 
has positive Hausdorff dimension, i.e. 
$$\Gamma_H^{ord}:=\{ y:~ \dim_{H}(L(y)) >0 \}.$$ 
Then $\Gamma_{H}^{ord}$ has 
full Hausdorff dimension, i.e.
\beql{183}
\dim_{H}( \Gamma_H^{ord})= 1.
\eeq
 It also has 
Lebesgue measure zero.
\end{theorem} 

\paragraph{Proof.}
The Hausdorff dimension $1$ property is shown  
by the author and Maddock \cite[Theorem 1.5]{LM10b}). The Lebesgue measure
zero property follows from Theorem~\ref{th94}(1).
$~~~\Box$.\medskip

There is potentially  a multifractal formalism connected to the
Hausdoff dimensions of the sets
$\Gamma_H^{ord}(\alpha):=\{ y:~ \dim_{H}(L(y)) >\alpha \},$ 
see \cite[Sec. 1.2]{LM10b}, and for multifractal formalism
see Jaffard  \cite{Ja97a}, \cite{Ja97b}.

Theorem \ref{th94} asserts that almost  all level sets in the ordinate sense are finite. 
The finite level sets have an intricate structure, which is analyzed by
Allaart \cite{A11b}. Finite level sets must have even cardinality, and all
even values occur.


\section{Local Level Sets of the Takagi Function} \label{sec13}

The  notion of {\em local level set} of the Takagi function
was recently introduced by the author and Maddock \cite{LM10a}.
These sets  are closed subsets
of level sets that are directly constructible from the binary expansions of any one of their
members. 

Local level sets are defined
by equivalence relation on elements $x \in [0,1]$ based on properties
of their binary expansion: $x= .b_1 b_2b_3\cdots$.
Recall that in Section \ref{sec2} we defined the deficient digit function
$$
\dN_j(x):= j -2(b_1+b_2+ \cdots + b_j).
$$
The quantity $\dN_j(x)$ counts the excess of binary digits $b_k=0$ over those with $b_k=1$
in the first $j$ digits. We then associate  to any $x$  the sequence  of  ``breakpoints" $j$
at which tie-values $\dN_j(x)=0$ occur, setting
$$
Z(x) : = \{ c_m:~~\dN_{c_m}(x)=0\}.
$$
where we define $c_0=c_0(x) = 0$ and set $c_0(x)< c_1(x)< c_2(x) < ...$. This sequence of tie-values may
be finite or infinite, and if it is finite, ending in $c_{n}(x)$, we make the
convention to adjoin a final ``breakpoint"  $c_{n+1}(x)= +\infty$.
 We call a ``block" a set of digits between two consecutive tie-values, 
$$
B_k(x) := \{ b_j: ~c_k(x)  < j \le c_{k+1}(x)\}.
$$
Two blocks are  called equivalent, 
written  $B_k(x) \sim B_{k'}(x')$, if their  endpoints agree
($c_k(x)= c_{k'}(x')$ and $c_{k+1}(x) = c_{k'+1}(x')$)
and either $B_k(x) = B_{k'}(x')$ or $B_k(x) = \bar{B}_{k'}(x')$, where the bar operation flips
all the digits in the block, i.e. 
$$b_j \mapsto b_j^{'}:= 1- b_j,~~~~~~~ c_k < j \le c_{k+1}.$$
Finally, we define the equivalence relation on binary expansions $x \sim x'$ to mean that
$Z(x) \equiv Z(x')$, and furthermore every block $B_k(x) \sim B_k(x')$ for $k \ge 0$.
We define the  {\em local level set} $L_x^{loc}$ associated to $x$
 to be the set of equivalent points, 
$$
L_x^{loc} := \{ x': ~~x' \sim x\}.
$$
It is easy to show using Takagi's formula (Theorem ~\ref{le21})
 that the relation $x\sim x'$ implies that $\tau(x)= \tau(x')$ so that $x$ and $x'$
are in the same level set of the Takagi function.

Each local level set $L_x^{loc}$  is a closed set, 
and is either a finite set if $Z(x)$ is finite, or else  is a perfect
totally disconnected set (Cantor set) if $Z(x)$ is infinite.
The case of dyadic rationals $x= \frac{k}{2^n}$ is exceptional, since they
have two binary expansions, and we remove this ambiguity by taking the binary
expansion for $x$ that ends in zeros. \\

The definition implies  that each level set $L(y)$ partitions into  a disjoint union of local
level sets $L_x^{loc}$. A priori this union may be finite, countable or uncountable.
Finite and countable examples are given by the author and Maddock \cite{LM10b}. For the  uncountable case, see
Theorem ~\ref{th125} below. Note that all countably infinite level sets  necessarily are a 
countable union of finite local level sets. The only level sets  currently
known to be countably infinite are certain dyadic rational levels, including $x= \frac{1}{2}.$

The paper \cite{LM10a} characterizes the size of a ``random" local level set sampled by randomly drawing
an abscissa value $x$, as follows; this result immediately implies Theorem ~\ref{th61}.

%
%
\begin{theorem}\label{th101}
For a full Lebesgue measure set of abscissa points $x \in [0,1]$ the local level set
$L_x^{loc}$ is a Cantor set of Hausdorff dimension $0$.
\end{theorem}

\paragraph{Proof.}
This is shown in  \cite[Theorem 1.4]{LM10a}.
$~~~\Box$. \smallskip

In order to analyze local  level sets in the ordinate space $0 \le y \le \frac{2}{3}$,
we label each local level set by its leftmost endpoint. 
We define the {\em deficient digit set} $\ddL$ by the condition
\beql{1221}
\ddL := \{ x \in [0,1]:~x=0.b_1b_2b_3... \mbox{such~that}~~ \dN_j(x) \ge 0, ~j=1, 2, 3, ...\}
\eeq
This set turns out to be quite important  for understanding the Takagi function.

%
%

\begin{theorem}\label{th102}
(1) The  deficient digit set 
$\ddL$ is the set of leftmost endpoints of all local level sets.

(2) The set $\ddL$ is a closed, perfect set (Cantor set). 
It has Lebesgue measure $0$.

(3) The set $\ddL$ has  Hausdorff dimension $1$.
\end{theorem}

\paragraph{Proof.}
(1) and (2) are shown in \cite[Theorem 4.6]{LM10a}.

(3) This was shown in \cite[Theorem 6.1]{LM10b}.
$~~~\Box$ \medskip

We show that the Takagi function restricted to the set
$\ddL$ is quite nicely behaved, as given in the
following result.

%
%

\begin{theorem} \label{th123} 
The function  $\tauS(x)$ defined by $\tauS(x)= \tau(x) + x$ for
 $x \in \ddL$  
is a nondecreasing function on $\ddL$. Define its extension to
all $x \in [0,1]$ by 
$$
\tauS(x) := \sup\{ \tauS(x_1):  x_1 \le x ~~~\mbox{with}  ~~ x_1 \in \ddL\}.
$$
 Then the  function $\tauS(x)$  is  a monotone singular function.  That is, it is
 a nondecreasing continuous function  
 having  
 $\tauS(0)=0, \tauS(1)=1$, which has derivative zero at (Lebesgue) almost all points of $[0,1]$. 
  The closure of the set of points of increase of $\tauS(x)$ is the deficient digit set $\ddL$.
\end{theorem}

\paragraph{Proof.} This is shown in \cite[Theorem 1.5]{LM10a}. $~~~\Box$.\smallskip

%
%

\begin{figure}[h]
  \begin{center}                                                               
   $ \begin{array}{c}
           \includegraphics[width=2.0in]{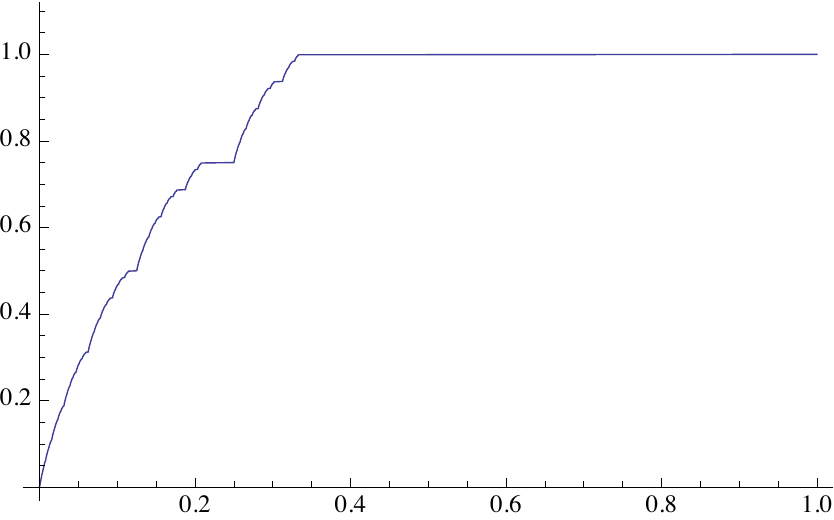} 
                    \end{array}$
  \end{center}
  \caption{Graph of 
    Takagi singular function 
 $\tauS(x)$}
  \label{fig121}
\end{figure}

We call the function $\tauS(x): [0,1] \to [0,1]$ the {\em Takagi singular function};
it  is pictured in Figure \ref{fig121}.
Using the functional equation
 $\tau(\frac{1}{2}x) = \frac{1}{2} (\tau(x) + x)$
given in Theorem \ref{th22}(1), we deduce that the function $\frac{1}{2} \tauS(x)$ agrees 
with $\tau(x)$
on the set $\frac{1}{2} \ddL$. This  shows that the Takagi function is
strictly increasing when restricted to the domain $\frac{1}{2} \ddL$.\smallskip

Associated to
the Takagi singular function is a nonnegative Radon measure $d\muT$, 
which we call  the {\em Takagi singular measure},
such that 
\beql{129a}
\tauS(x) = \int_{0}^x d\muT.
\eeq
It is singular with respect to Lebesgue measure and
defines an interesting probability measure on $[0,1]$.
The Takagi singular  measure is not translation-invariant,
but it has  certain self-similarity properties under dyadic
rescalings, compatible with the functional equations of
the Takagi function. 
These are useful in explicitly computing
the measure of various interesting subsets of $\ddL$, see \cite{LM10b}. 
One may compare its properties with those 
of the Cantor function, see Dovghoshey et al \cite[Sect. 5]{DMRV06}.
A major difference is that its support $\ddL$ has full Hausdorff dimension,
while the Cantor set has Hausdorff dimension $\log_3 (2) \approx 0.63092.$\smallskip

By  an application of 
the co-area formula of geometric measure theory for functions of bounded variation
(in the version in Leoni \cite[Theorem 7.2 and Theorem 13.25]{Le09}) 
 to a relative of the Takagi singular function, the author with Maddock \cite{LM10a}  determined
the average number of local level sets on a random level.

%
%
\begin{theorem}\label{th124}
For a full Lebesgue measure set of ordinate points $y \in [0, \frac{2}{3}]$ the
number $N^{loc}(y)$ of local level sets at level $y$ is finite. Furthermore
$$
\int_{0}^{\frac{2}{3}}  N^{loc}(y) dy = 1.
$$
That is,  the expected number of local level sets 
at a uniformly drawn random level in $[0, \frac{2}{3}]$ is exactly $\frac{3}{2}$.
\end{theorem}

\paragraph{Proof.} This is shown in \cite[Theorem 6.3]{LM10a}. 
A simplified derivation is given in \cite{A11},
that avoids the co-area formula $~~~\Box$\smallskip

One can easily derive  both parts of Theorem~\ref{th94}  from this result. 
This result  fails to give any information 
about the multiplicity of local level sets on those levels having an uncountable level
set, because the set of such levels $y$ has Lebesgue measure $0$.\smallskip

Recently Allaart \cite{A11} established the following result
about the multiplicities  of local level sets in a level set.

%
%
\begin{theorem}\label{th125} {\rm (Allaart 2011)}
There exist levels $y$ such that
the level set $L(y)$ contains  uncountably many distinct local level sets.
The set of such levels is dense in $[0, \frac{2}{3}].$
\end{theorem}

Allaart (\cite{A11} , \cite{A11b}) has obtained further information on the structure
of level sets and local level sets, classifying the values giving different types of level sets
in terms of the Borel hierarchy of descriptive set theory. Some of these sets
are complicated enough that they apparently live in the third level of the 
Borel hierarchy.


%
%
%

\section{Level Sets at Rational Levels}\label{sec14}

In 2005 Knuth~\cite[Problem 83, p. 32]{Kn05} raised the question of determining which
rational levels $y= \frac{r}{s} $ have level sets have $L(\frac{r}{s})$
that are  uncountable. 

The author with  Maddock \cite{LM10a} answered  the much easier question 
of determining  when certain rational numbers $x$ give uncountable local
level sets $L_x^{loc} \subset L(\tau(x))$. 

%
%

\begin{theorem}~\label{th131}
(1) A rational number $x=\frac{p}{q} \in [0,1]$ has an uncountable local level
set $L_{\frac{p}{q}}^{loc}$ if and only if  its binary expansion
has a pre-periodic part with an equal number of zeros and ones, and if
also its purely periodic part has an equal number of zeros and ones.

(2) If a rational $x$ has $L_x^{loc}$
uncountable, then $L_x^{loc}$ contains a countably infinite
set of rational numbers. 
\end{theorem}


This criterion implies that  dyadic rationals $x= \frac{m}{2^n}$ must have finite local level sets. 
For such values  $y=\tau(\frac{k}{2^l})$ is also a dyadic rational. Concerning
dyadic rational levels, 
Allaart \cite{A11b} obtains the following much stronger result.

%
%

\begin{theorem}~\label{th112}{\rm (Allaart  2011)}
Let $y= \frac{k}{2^m}$ be a dyadic rational with $0 \le y \le \frac{2}{3}$.
Then the level set $L(y)$ is either finite or countably infinite. 
\end{theorem}

The examples   $y=0$ and $y= \frac{1}{2}$ given above show
that both alternatives in this result  occur. 
Furthermore all elements on the given dyadic rational level are rational. 
However if  the set is countably infinite  then these elements  need not all be dyadic rationals. For
example $\tau(\frac{1}{6})= \frac{1}{2},$ cf. \cite[Theorem 7.1]{LM10a}.\smallskip

%
%
%

\section{Open Problems}\label{sec15}

There remain many open questions 
about the Takagi function.  We mention a few of them.

\begin{enumerate}
\item
Determine analytic and other
properties of the conjugate Takagi
function $Y(x)$ defined in Theorem~\ref{th32a}.
\item
Consider the set of all $x$ such that 
the abscissa level set  $L(\tau(x))$ has
Hausdorff dimension zero. Does this set have
full Lebesgue measure in $[0,1]$?

\item
Is there a continuous function on the interval for which the Takagi singular measure
is a ``natural"  invariant probability measure?
\item
Determine the dimension spectrum function
$$
f_{\tau}(\alpha) := \dim_H( \{ y:  ~~\dim_{H}( L(y)) \ge \alpha\}.
$$
(Some bounds on  $f_{\tau}(\alpha)$ are given in \cite{LM10b}.)
\item
Knuth's problem (\cite[7.2.1.3, Prob. 83]{Kn05}): Find a direct characterization of  the set of  rational $y$ which have 
an uncountable  level set. 
\end{enumerate}

\paragraph{Acknowledgments.} 
 D. E. Knuth originally brought  problems on the Takagi
function to my attention via postcard.  I thank S.T. Kuroda
for  helpful remarks on Takagi's original work,
and J.-P. Allouche for some corrections. 
I thank Pieter Allaart and Zachary Maddock for  helpful comments and for
supplying additional  references. I  gratefully thank the reviewer
for detailed comments and corrections. 

%
%
%
%

\noindent Jeffrey C. Lagarias \\
Dept. of Mathematics\\
The University of Michigan \\
Ann Arbor, MI 48109-1043\\
\noindent email: {\tt lagarias@umich.edu}\\

\end{document}